\documentclass[ijds,nonblindrev]{informs4}

\OneAndAHalfSpacedXI

\usepackage{amsmath,amssymb,amsfonts,easyReview}
\usepackage{appendix}
\usepackage{hyperref}
\usepackage{booktabs}
\hypersetup{
    hidelinks,
    colorlinks = true,
    linkcolor = blue,
    citecolor = blue,
    urlcolor = green
}
\usepackage[ruled,linesnumbered]{algorithm2e}

\usepackage{natbib}
 \bibpunct[, ]{(}{)}{,}{a}{}{,}%
 \def\bibfont{\small}%

\usepackage{rotating}
\usepackage{fancyvrb}
\TheoremsNumberedThrough     
\ECRepeatTheorems
\JOURNAL{Any INFORMS Journal}

\EquationsNumberedThrough    

\MANUSCRIPTNO{IJDS-0001-1922.65}

\begin{document}


\RUNTITLE{Zhaoen Li et al.}

\TITLE{Achieving Robust Data-driven Contextual
Decision Making in a Data Augmentation Way}

\ARTICLEAUTHORS{%
\AUTHOR{Zhaoen Li}

\AFF{Department of Industrial Engineering, Tsinghua University, Beijing 100084, China, \EMAIL{lze21@mails.tsinghua.edu.cn}}
\AUTHOR{Maoqi Liu}

\AFF{School of Management, Shandong University, Jinan 250100, China, \EMAIL{liumq@sdu.edu.cn}}
\AUTHOR{Zhi-Hai Zhang$^{\ast}$}

\AFF{Department of Industrial Engineering, Tsinghua University, Beijing 100084, China, \EMAIL{zhzhang@tsinghua.edu.cn}}

}

\ABSTRACT{%
    This paper focuses on the contextual optimization problem where a decision is subject to some uncertain parameters and covariates that have some predictive power on those parameters are available before the decision is made. More specifically, we focus on solving the Wasserstein-distance-based distributionally robust optimization (DRO) model for the problem, which maximizes the worst-case expected objective over an uncertainty set including all distributions closed enough to a nominal distribution with respect to the Wasserstein distance. We develop a stochastic gradient descent algorithm based on the idea of data augmentation to solve the model efficiently. The algorithm iteratively a) does a bootstrapping sample from the nominal distribution; b) perturbs the adversarially and c) updates decisions. Accordingly, the computational time of the algorithm is only determined by the number of iterations and the complexity of computing the gradient of a single sample. Except for efficiently solving the model, the algorithm provide additional advantages that the proposed algorithm can cope with any nominal distributions and therefore is extendable to solve the problem in an online setting. We also prove that the algorithm converges to the optimal solution of the DRO model at a rate of a $O(1/\sqrt{T})$, where $T$ is the number of iterations of bootstrapping. Consequently, the performance guarantee of the algorithm is that of the DRO model plus $O(1/\sqrt{T})$. Through extensive numerical experiments, we demonstrate the superior performance of the proposed algorithm to several benchmarks. 
}


\maketitle
\section{Introduction}\label{sec1}
Contextual optimization stands out as a prominent area bridging operations research and machine learning, characterized by decision-making scenarios with uncertain parameters whose distribution is initially unknown but can be partially inferred through observed covariates prior to the decision \citep{sadana2023survey}. This concept finds broad applicability across various domains, including inventory control \citep{ban2019big,bertsimas2022data,zhang2023optimal}, revenue management \citep{cohen2022price,liu2023active,perakis2023robust}, and portfolio optimization \citep{wang2021distributionally,elmachtoub2022smart}. In such contextual decision-making problems, the decision maker aims to address the conditional stochastic programming:
\begin{align}
	\mathop{\inf} \limits_{z \in \mathcal{Z}} \mathbb E_{\boldsymbol{Y}}[c(z;\boldsymbol{Y})|\boldsymbol{X}=\boldsymbol{x}],\forall \boldsymbol{x}\in\mathcal{X}. \label{e1}
\end{align}
Here, $c(\cdot)$ represents the cost function influenced by decisions denoted as $z$ and random parameters denoted as $\boldsymbol{Y} \in \mathcal{Y}$, which are unobservable during decision-making. The vector $\boldsymbol{X}$ signifies covariates that impact the distribution of $\boldsymbol{Y}$ and are known during decision-making processes.

Traditionally, contextual optimization follows a two-step process: first, establishing a predictive model for uncertain parameters, and then integrating this prediction into subsequent optimization procedures. However, a significant drawback of this method is its inability to consider the influence of prediction errors on subsequent optimizations. Consequently, even minor prediction errors can result in substantial sub-optimality during decision-making \citep{qi2021integrated,elmachtoub2022smart,sadana2023survey}. To address this issue, extensive efforts have focused on end-to-end approaches based on historical observations of covariates and uncertain parameters. These approaches aim to directly link covariates to optimal decision outcomes, thereby reducing the impact of prediction errors. Among these approaches, a commonly employed strategy is training a decision policy using observed covariates, which can be formally expressed as follows: 
\begin{align}
	\mathop{\inf} \limits_{\theta \in \Theta} \mathbb E_{\boldsymbol{X},\boldsymbol{Y}}[c(f(\theta;\boldsymbol{X});\boldsymbol{Y})], \label{e2}
\end{align}
where $f(\theta;\boldsymbol{X})$ represents the function mapping from the feature space to the decision space with coefficients $\theta$ in the space $\Theta$. Using the interchangeability principle \citep{zhang2023optimal}, we have that 
\begin{equation}
    \mathbb E_{\boldsymbol{X}}\Big [\mathop{\inf} \limits_{z \in \mathcal{Z}} \mathbb E_{\boldsymbol{Y}}[c(z;\boldsymbol{Y})|\boldsymbol{X}] \Big] = \mathop{\inf} \limits_{\theta \in \Theta} \mathbb E_{\boldsymbol{X},\boldsymbol{Y}}[c(f(\theta;\boldsymbol{X});\boldsymbol{Y})]. \notag
\end{equation}
That is to say, by taking the expectation of \eqref{e1} over the marginal distribution of the feature ${\boldsymbol{X}}$, the optimal policy $f(\theta^{\star};\boldsymbol{X})$ obtained in \eqref{e2} minimizes the marginalized expected cost of \eqref{e1}. The existing literature has presented numerous analytical and numerical experiments showcasing the substantial performance enhancements achieved by problem \eqref{e2} in contrast to the traditional predict-then-optimize approach \citep{ban2019big,bertsimas2022data,bertsimas2023dynamic,rychener2023end}.

However, formulating and solve model \eqref{e2} pose twofold challenges. The first is the imperfect knowledge on the joint distribution. The decision maker is rare to, if not impossible, get access to the ground-truth joint distribution of $(\boldsymbol{X},\boldsymbol{Y})$. In contrast, only a noisy approximation, such as the empirical distribution established from a historical sample of $(\boldsymbol{X},\boldsymbol{Y})$, is available. The issues inspire the application of the distributionally robust optimization techniques in the contextual optimization. In general, the distributionally robust optimization (DRO) approach of problem \eqref{e2} can be formulated as
\begin{equation} 
	\mathop{\inf} \limits_{\theta \in \Theta} \mathop{\sup} \limits_{P \in \mathcal{P}} \quad \mathbb E_{P}[c(f(\theta;\boldsymbol{X});\boldsymbol{Y})], \label{e3}
\end{equation} 
where $\mathcal{P}$ represents the ambiguity set containing possible joint distributions of $(\boldsymbol{X},\boldsymbol{Y})$, denoted as $P$. A well-designed $\mathcal{P}$ can encompass the ground truth distribution with high probability, thereby offering a lower bound on problem \eqref{e2}.

The second challenge is the computational burden associated with solving problem \eqref{e2}. Compared to the traditional predict-then-optimize framework or non-contextual optimization, the computational complexity of contextual optimization is raised by the interplay of covariates, uncertainty parameters, and decision policies, compounded by the high dimensionality of the joint distribution of $(\boldsymbol{X},\boldsymbol{Y})$.

To address these challenges, we propose a novel data augmentation based stochastic gradient descent (DA-SGD) algorithm designed to solve the DRO model \eqref{e3} established with the Wasserstein uncertainty set, which encompasses distributions sufficiently close to a nominal distribution in terms of the Wasserstein distance. The algorithm iteratively updates the decision policy through three steps in each iteration. First, it conducts a bootstrapping draw from the nominal distribution to generate a synthetic sample.  Then, it adversarially perturbs this synthetic sample according to the objective of the DRO model with the coefficients determined in the previous iteration. Finally, it updates the coefficients of the decision policy based on the gradient computed at the perturbed synthetic sample. 
The contribution of our DA-SGD algorithm is summarized as follows. 
\begin{itemize}
	\item The proposed DA-SGD algorithm offers an efficient approach to solving Wasserstein-based DRO problems in contextual optimization \eqref{e3}. For arbitrary nominal distribution, the gap between optimal solution and DA-SGD is $O(\frac{1}{\sqrt{T}})+\nu$, where $T$ is the iterations of proposed algorithm and $\nu$ is the fixed gap between ground-truth distribution and nominal distribution. If nominal distribution is empirical distribution, the gap becomes $O(\frac{1}{\sqrt{T}}+\frac{1}{\sqrt{N}})$, where $N$ denotes the sample size. Considering the online learning scheme, the gap is $O(\frac{1}{\sqrt{T}})$, where $T$ represents both iterations and sample size.
 
    \item Our algorithm is rooted in the bootstrap method, allowing it to handle Wasserstein-based DRO models with arbitrary nominal distributions, instead of being restricted to specific distributions such as empirical or Elliptical distribution \citep{kuhn2019wasserstein, shafieezadeh2019regularization}. Based on this property, it is also well-suited for robust contextual optimization in an online scenario. Here, we can iteratively draw new samples from the ground-truth distribution, facilitating dynamic adaptation and continual improvement of the optimization process.
	
	\item In each iteration, our SGD algorithm perturbs the bootstrapped sample, aggregating which results in an augmented bootstrap distribution that is more adversarial to the decision policy. This comparison between the augmented and nominal distributions enables decision makers to gain deeper insights into the interplay among covariates, uncertain parameters, and decision policies.

    \item In numerical experiment, DA-SGD performs well in small data size than other benchmarks, and is least affected by the randomness of data generating. Compared with other DRO algorithm, it saves a lot of solving time to obtain the optimal solution.
	
\end{itemize}

The remainder of the paper is organized as follows: Section \ref{sec2} presents a comprehensive review of the related literature. In Section \ref{sec3}, we present the framework of DA-SGD algorithm and its online learning version. In Section \ref{sec4}, we provide the performance guarantee of new proposed algorithm, and show the finite sample guarantee and online regret in special cases. In Section \ref{sec5}, we conduct the experiment and report numerical results. Section \ref{sec6} draws conclusions and outlines future research directions. The proof of the propositions, lemmas and theorems in this paper can be found in Appendix.

\section{Literature Review}\label{sec2}
Contextual optimization (CO) has garnered extensive attention in previous literature, which can be categorized into three basic formats: decision rule optimization \citep{ban2019big,bertsimas2019optimal,zhang2023optimal}, sequential learning and optimization \citep{ferreira2016analytics,bertsimas2020predictive,bertsimas2022bootstrap}, and integrated learning and optimization \citep{elmachtoub2022smart,qi2021integrated,gong2023training}. In this paper, our focus is on decision rule optimization. For an in-depth review of the other two schemes, we refer readers to \cite{sadana2023survey}.

The decision rule optimization approach for contextual optimization was popularized by \cite{ban2019big} where a CO counterpart of newsvendor problem is studied. Through extensive analytical and numerical results, they showcased the benefits of incorporating covariate information into the problem. The linear decision rule, where the decision is a linear combination of the covariates, and its extensions have received extensive attention in previous literature.  \cite{bazier2020generalization} studied contextual portfolio optimization based on the linear decision rule with regularization and derived a performance bound for the problem. For general contextual optimization problems, \cite{bertsimas2022data} investigated linear policies in reproducing kernel Hilbert space (RKHS), and demonstrated the asymptotic optimality of RKHS-based policies. \cite{notz2022prescriptive} applied RKHS-based policies to the flexible capacity allocation problem and provided performance guarantees. \cite{zhang2023data} focused on piecewise-affine decision rules and proved the performance guarantee of the proposed model.

Other researchers have explored machine learning approaches as proxies for decision rules. \cite{keshavarz2022interpretable} trained a tree-based machine learning model for decision policies in a contextual newsvendor problem. \cite{bertsimas2019optimal} proposed a prescriptive tree model for general CO problems. \cite{huber2019data} utilized gradient boost decision trees and neural networks in CO of newsvendor problems. Similarly, \cite{oroojlooyjadid2020applying} applied deep neural networks to the problem. \cite{zhang2017assessing} modified the loss function of a deep neural network from quadratic to $L_1$-norm and demonstrated improved performance. \cite{zhang2021universal} proposed a contextual portfolio optimization model based on deep neural networks and discussed adjusting the output layer to incorporate basic constraints. To ensure feasibility in an economic dispatch problem, \cite{chen2023end} added a repaired layer to a deep neural network. \cite{ciocan2022interpretable} learned action policies for optimal stopping problems using decision trees. 

These approaches relied on an empirical risk minimization approach and often overlooked data uncertainty, leading to sub-optimal policies. Therefore, it is crucial to incorporate robustness against data uncertainty. Due to its natural relationship with Empirical Risk Minimization (ERM), the Distributionally Robust Optimization (DRO) approach based on statistical distance is applied to contextual optimization (CO). The key idea behind establishing the statistical-distance-based DRO model is to consider all distributions that are close enough to a nominal one with respect to some statistical distance.\cite{yang2022decision} proposed a DRO counterpart of CO based on the causal transport distance and proved its performance guarantee. \cite{zhang2023optimal} introduced a robust policy optimization framework based on the Wasserstein metric to obtain robust policies for the contextual newsvendor problem and used a method they called the Shapley extension to achieve the global optimal solution. \cite{rychener2023end} proposed a stochastic gradient descent (SGD) algorithm that is applied to solve the DRO model with $\phi$-divergence.

Our work distinguishes itself from these DRO studies in two significant ways. Firstly, the above models relied on the nominal distribution to be the empirical distribution to maintain tractability. Both \cite{yang2022decision} and \cite{zhang2023optimal} relied on the empirical nominal distribution to get the tracable reformulation of their proposed model. In contrast, our approach leverages bootstrapping techniques to accommodate any nominal distribution.  Secondly, we provide an SGD algorithm with desirable asymptotic optimality and finite sample guarantee. The scale of the dual problem is proportionate to the size of the data set used to obtain the empirical distribution. Although \cite{rychener2023end} also proposed an SGD algorithm for the DRO model, their algorithm's performance is not guaranteed for the Wasserstein DRO, and it requires solving the DRO model in each iteration. In contrast, our algorithm updates the solution using one sample in each iteration and is proven to be asymptotically optimal, converging at a rate of $O(\sqrt{N})$, where $N$ is the number of available samples used in bootstrap.

\section{DA-SGD Framework for DRO}\label{sec3}
In this section, we first specify the uncertainty set based on Wasserstein distance and provide some preliminary of the corresponding DRO model in Section \ref{sec3.1}. And then, we formally articulate the proposed SGD algorithm to solve Problem \eqref{e3} and demonstrate how to modify it to an online learning problems in Section \ref{sec3.2}. 

\subsection{Preliminary on Wasserstein-distance-based Distributionally Robust Formulation} \label{sec3.1}
The core idea of Wasserstein-distance-based DRO model is to maximize the worst-case expectation over all distribution closed to a nominal one where the distance between distributions is measured by the Wasserstein distance, and the distribution set containing all the distribution above is called ambiguity set. Formally, Wasserstein distance is defined as 
\begin{definition}[Wasserstein metric, \cite{mohajerin2018data}]\label{WasMetric}
	The Wasserstein distance $W_d:\mathcal{M}\left(\Xi\right)\times
	\mathcal{M}\left(\Xi\right) \rightarrow \mathbb{R_{+}}$ is defined as,
	\begin{equation}
		W_d\left(P,P_0\right)=\inf\limits_{\Pi}\left\{ \int_{\Xi^2} d(\xi,\xi_0) \Pi\left(d\boldsymbol \xi,d\boldsymbol \xi_0 \right): \left.
		\begin{array}{l}
			\textit{$\Pi$ \rm{is a joint distribution distribution of} $\boldsymbol \xi$ \rm{and} $\boldsymbol \xi_0$}\\
			\textit{ \rm{with marginals} $P$ \rm{and} $P_0$, \rm{respectively}}
		\end{array}
		\right.\right\}, \nonumber
	\end{equation}
	where $\Xi^2$ is the support of $\boldsymbol \xi$ and $\boldsymbol \xi_0$; $\mathcal{M}\left(\Xi\right)$ is a set of all distributions that supported by $\Xi$; and 	$d(\xi,\xi')$ is some distance metric. In this study, $\xi=(x,y)$ is the sample point, and we consider a norm-based mass transportation cost in Wasserstein distance, which is: 
	\begin{equation}
		d(\xi,\xi')=\|x-x'\|_2 + \kappa\|y-y'\|_2, \notag
	\end{equation}
	for some $\kappa>0$. $\kappa$ controls the costs of moving probability mass along the label space. For $\kappa=\infty$, all distribution in the ambiguity set are obtained only along the feature space \citep{shafieezadeh2019regularization}. In this study, $\kappa$ is set as $\infty$.
\end{definition} 

Let $Q$ denotes the nominal distribution and then Wasserstein distance based ambiguity set $\mathcal{P}$ has the form as:
\begin{equation}
	\mathcal{P} = \{P:W_d(P,Q)\leq\rho\}, \notag
\end{equation}
where $\rho$ is an upper bound of Wasserstein distance between a possible distribution and the nominal one. However, even with these definitions, Problem \eqref{e3} is still hard to solve, and one assumption should be made then.

\begin{assumption} \label{asp1}
    The objective function and distribution are all defined on a compact support space $\Xi$ with the diameter $D$, where
    \begin{equation}
        D=\mathop{\sup}\{d(\xi_i,\xi_j):\xi_i,\xi_j\in\Xi\}. \notag
    \end{equation}
\end{assumption}

As shown in \cite{gao2023distributionally}, the Wasserstein-distance-based DRO model has a tractable reformulation by its strong duality. We adopt the result to the CO context and get the tractable reformulation of Problem \eqref{e3} as follows when Assumption \ref{asp1} holds:
\begin{equation}
	\mathop{\inf} \limits_{\theta \in \Theta} \mathop{\sup} \limits_{P \in \mathcal{P}} \mathbb	 E_{P}[c(f(\theta;\boldsymbol{X});\boldsymbol{Y})]=\mathop{\inf} \limits_{\theta \in \Theta, \gamma>0} \Big\{\mathbb{E} _{Q}[\mathop{\sup}\limits_{\xi' \in \Xi} \{c(f(\theta;X');Y')-\gamma d(\xi',\xi)+\gamma\rho\}]\Big\},  \label{e4}
\end{equation}
and to simplify the formulation above,  define
\begin{equation} \label{e5}
    h(\theta,\gamma;\xi,\xi') = c(f(\theta;X');Y')-\gamma d(\xi',\xi)+\gamma\rho,
\end{equation}
\begin{equation}\label{e6}
    h^{\star}(\theta,\gamma;\xi) = \mathop{\sup}\limits_{\xi' \in \Xi} h(\theta,\gamma;\xi,\xi'),
\end{equation}
\begin{equation} \label{e7}
    H(\theta,\gamma) = \mathbb{E} _{Q}[h^{\star}(\theta,\gamma;\xi)].
\end{equation}

Without further assumption on $Q$, solving problem \eqref{e4} encounters complicated integrals. Typically, the nominal distribution is assumed to be empirical or Gaussian distributions to further simplify the calculation \eqref{e4}.  In the following section, we propose an algorithm based on stochastic gradient descent to solve problem \eqref{e4} without any limitations on the nominal distribution.

\subsection{Data Augmentation based Stochastic Gradient Descent Algorithm} \label{sec3.2}
In Section \ref{sec3.1}, problem \eqref{e3} is reformulated as a min-max problems \eqref{e5}. The inner-maximization has an intuitive explanation that for a given $\theta$, $\gamma$ and the nominal distribution $Q$, the mass at $\xi\in\Xi$ is perturbed to $\xi'\in\Xi$ to get the worst-case sample. This perturbation is in the same spirit of the adversarial training based on data augmentation in the machine learning field, which generates adversarial synthetic data based on the bootstrap techniques and improve the generalization level of prediction model \citep{goodfellow2014explaining,engstrom2019adversarial}. We take advantages of the idea to develop an iterative stochastic gradient descent algorithm to solve the DRO counterpart of CO in this section.The detail of DA-SGD algorithm is presented in Algorithm \ref{algorithm1}. 
\begin{algorithm}[htbp] 
	\KwData{nominal distribution $Q$, radius $\rho$ and threshold $\epsilon$}
	\KwResult{learned model weights $\theta$}
	initialization model weights $\theta=\theta_0$ and $\gamma=\gamma_0$\ randomly;\\
	\For{$t=1,...,T$}{
		Sample $\xi_t$ according to $Q$\;
		$\hat{\xi}_t^{\star} \leftarrow \xi_t$; \\
		$k\leftarrow0$;\\
		\While{$\nabla_{\xi}h(\theta,\gamma;\xi_t,\hat{\xi}_t^{\star})>\epsilon$ and $k\leq K$}{
			$\hat{\xi}_t^{\star}\leftarrow\hat{\xi}_t^{\star}+\eta\nabla_{\xi}h(\theta,\gamma;\xi_t,\hat{\xi}_t^{\star})$\\
			$k\leftarrow k+1$\\
		}
		$\theta\leftarrow\theta-\alpha\nabla_{\theta}c(f(\theta;\hat{x}_t^{\star});\hat{y}_t^{\star})$\\
		$\gamma\leftarrow\gamma-\beta\bigg[-d(\hat{\xi}_t^{\star},\xi_t)+\rho\bigg]$ \\
	}
	\caption{DA-SGD algorithm}
	\label{algorithm1}
\end{algorithm}

To apply this algorithm, except for the nominal distribution $Q$ and the radius of Wasserstein uncertainty set $\rho$, some parameters are also needed. $\epsilon$ is a threshold to bound the gap between $\hat{\xi}_t^{\star}$ and real optimal point $\xi_t^{\star}$. In augmentation part, worst-case sample is generated based on the inner maximization problem. Using 2-norm distance in it, this problem is time consuming when the dimensions of dataset increasing. Thus, we also apply a gradient descent method for it, which can save a lot of time and limit the time with the help of $K$. To avoid the situation when trapped in one point, $K$ is set to terminate augmentation procedure in time. $\eta$ is the learning rate in augmentation part, and $\alpha$ and $\beta$ are the learning rate for gradient descent, which are set in Theorem \ref{thm1}. The core of Algorithm \ref{algorithm1} is separated into three key components. The first one is release of the decision maker to know the nominal distribution in advance, and then sample point $\hat{\xi}$ from the nominal distribution following bootstrap techniques. The second is the augmentation of the sampled point according to the maximization $\mathop{\sup}\limits_{\xi' \in \Xi} h(\theta,\gamma;\xi,\xi')$. Finally, we use augmented data point to update $\theta$ and $\gamma$ via gradient descent with the learning rate $\alpha$ and $\beta$.

The bootstrap approach in the algorithm free it from the specific form of the nominal distribution $Q$, which leverage the algorithm twofold advantages. 
\begin{itemize}
    \item Firstly, our algorithm can reduce computational burden. In each iteration, only one sample participate in the inner maximization problem and gradient descent, which saves lots of time for solving initial problem and makes it more possible to find the global optimal point. The computational time increases linearly in the number of bootstrap $T$. In the next section, we show that the algorithm converge to the optimal solution to the DRO \eqref{e3} with a rage of $O(\frac{1}{\sqrt{T}})$, which leverage a desirable efficiency of the algorithm to solve problem \eqref{e3}. 
    \item Secondly, the algorithm do not requires us to know the full information of the nominal distribution in advance of the implementation as long as we can draw i.i.d samples from the nominal distribution. This is advance of taking the ground truth distribution as the nominal one and each iteration as drawing a new sample from the ground truth distribution. As such, the algorithm is applicable to an online context where the data stream keep coming into the system and the decision is continuously done as each data point is observed. The context is quite common in CO such as the recommendation system, economic dispatching and many other fields.
\end{itemize}

In each iteration of Algorithm \ref{algorithm1}, we first do the bootstrap to overcome small data problems and sample one point from the dataset. Then the sampled point is augmented to solve the inner maximization problem. Finally, we use augmented data point to update $\theta$ and $\gamma$ via gradient descent. In terms of computational complexity, since algorithm \ref{algorithm1} does not need to solve optimization problems directly, it can greatly reduce the solving time. And the computational time of our algorithm is $O(T\cdot K)$, which is only related to the data augmentation iterations and training iterations, and is not affected by the size of dataset. That is to say, much solving time can be saved by slightly reducing the quality of the solution. 

\begin{remark}[Applying Algorithm \ref{algorithm1} to an online learning context]
    In Algorithm \ref{algorithm1}, if the bootstrap sampling step is replaced by data coming which is described in online learning method, then Algorithm \ref{algorithm1} becomes an online learning version algorithm.
    Because we do not need the specific nominal distribution when implementing algorithm \ref{algorithm1}, it enable us to apply the online learning context where the ground truth distribution is unknown or even changeable to the decision maker, while i.i.d samples can be sequentially drawn from the distribution as time goes by. In this situation, we can consider $t$ as the time when $t-$th sample arrives, and time $T$ is the total time.
 \end{remark}

\section{Performance Guarantee} \label{sec4}
In this section, we will show the theoretical guarantee of the DA-SGD algorithm. That is, we aim to characterize the performance gap between the solutions to the proposed DA-SGD algorithm and the optimal one to problem \eqref{e2} the oracle distribution denoted as $P_{true}$. Let $\theta_T,\theta_{R},\theta^{\star}$ denote the oracle optimal solution under $P_{true}$, that to the proposed the proposed algorithm after $T$ iterations, and that to the DRO model \eqref{e3}. Let $C(\theta)=\mathbb{E}_{P_{true}}[c(f(\theta;X);Y)]$, and $C_{DRO}(\theta;\rho)$ be the oracle objective of a given solution $\theta$ and the DRO model with given $\theta$ and robust parameter $\rho$, respectively;  Specifically, we will investigate the following gap $r$ in this section:
\begin{equation}
    r = C(\theta_T) - C(\theta^{\star})
\end{equation}
In the following, we will analyze the gap from two aspects:
\begin{itemize}
    \item [a)] the gap between the solutions to the SGD algorithm and the DRO model \eqref{e3}, which is defined as
    $$r_{algorithm}= C(\theta_{T})-C(\theta_{R}).$$
    \item [b)] the gap between the solution to the DRO model \eqref{e3} and the oracle optimal one.
    $$r_{DRO}=C(\theta_{R})-C(\theta^{\star})$$
\end{itemize}

Accordingly the gap $r$ can be reformulated into
\begin{equation}
    r = r_{algorithm} + r_{DRO}
\end{equation}

Since many literature has given the bound on $r_{DRO}$, the above division helps us to take advantages on the current results and clearly characterize the additional gap brought by the algorithm. 


In Section \ref{sec4.1}, we provide a upper bound on $r_{algorithm}$ as a function of $T$. In Section \ref{sec4.2}, we focus on a special case where the nominal distribution is the empirical distribution as an example of deriving $r_{DRO}$ and subsequently, provide a bound on $r$ in such a case. In Section \ref{sec4.3}, an upper bound of $r_{DRO}$ is given in online learning scheme.

\subsection{The Gap between Algorithm Solutions and DRO Model} \label{sec4.1}

In each iteration of the proposed algorithm, the worst sample-point $\xi^{\star}$ is determined by the values of both $\theta$ and $\gamma$ at that iteration, and therefore, can be described as a function $\xi^{\star}(\theta,\gamma)$ of $\theta$ and $\gamma$, which makes $h^{\star}(\theta,\gamma;\xi)=h(\theta,\gamma;\xi,\xi^{\star}(\theta,\gamma))$. In order to get the convergence about algorithm, we need first to bound the gradient of $h^{\star}(\theta,\gamma;\xi)$ and make one assumption about Lipschitz smoothness on objective function.

\begin{assumption} \label{asp2}
    Objective function $c(f(\theta;X);Y)$ satisfies Lipschitz smoothness conditions in both $\theta$ and $\xi$, where we have 
	\begin{equation}
		\begin{aligned}
			\|\nabla_{\theta}c(f(\theta;X);Y)-\nabla_{\theta}c(f(\theta';X);Y)\|_*\leq L_{\theta\theta}\|\theta-\theta'\|, \\
			\|\nabla_{\theta}c(f(\theta;X);Y)-\nabla_{\theta}c(f(\theta;X');Y')\|_*\leq L_{\theta\xi}\|\xi-\xi'\|, \\
			\|\nabla_{\xi}c(f(\theta;X);Y)-\nabla_{\xi}c(f(\theta';X);Y)\|_*\leq L_{\xi\theta}\|\theta-\theta'\|, \\
			\|\nabla_{\xi}c(f(\theta;X);Y)-\nabla_{\xi}c(f(\theta;X');Y')\|_*\leq L_{\xi\xi}\|\xi-\xi'\|,
		\end{aligned}
	\end{equation}
 where $\|\|_*$ denotes the dual norm, $L_{\theta\theta}$, $L_{\theta\xi}$, $L_{\xi\theta}$ and $L_{\xi\xi}$ are the Lipschitz constant for $c(f(\theta;X);Y)$ with respect to the combinations of $\theta$ and $\xi$.
\end{assumption}

This assumption is first proposed in \cite{sinha2017certifying} and applies to extensive settings in previous literature. For example, in linear programming problems, where $c(z;Y)$ is linear in $z$ and $Y$, the following examples are some cases that satisfy Assumption \ref{asp2}.

\begin{example} (Linear Policy)
    If linear policy is applied to predict the final decision, i.e.,
    $$f(\theta;X)=\theta^T g(X),$$
    Assumption \ref{asp2} holds when Assumption \ref{asp1} holds and $g(X)$ is a second-order differentiable function of $X$.
\end{example}

\begin{example} (Piece-wise Linear Policy)
    When considering piece-wise linear policy for decision-making with 
    $$f(\theta;X)=\max_{k\in[K]}(\theta^k)^{\top}X,$$ 
    and utilize a quadratic function to smooth the breakpoints, then Assumption \ref{asp2} holds. A detailed case of this example can be found in Section \ref{sec5}.
\end{example}

\begin{example} (Machine Learning)
    If Assumption \ref{asp1} holds, let $x'=x+\delta_0$ ($\delta_0$ always exists when $\Xi$ is compact) and $x'_i \neq 0$, then Sigmoid function
    \begin{equation}
        f(\theta;X) = \frac{1}{1+e^{-\theta^{\top}X}},
    \end{equation}
    and ELU function 
    \begin{equation}
        f(\theta;X) = \left\{
		\begin{aligned}
			& \theta^{\top}X &&if &&\theta^{\top}X>0\\
			&\alpha(e^{\theta^{\top}X}-1)&&if &&\theta^{\top}X\leq 0,\\
		\end{aligned}
		\right.
    \end{equation}
    satisfy Assumption \ref{asp2}.
\end{example}

Besides this, in order to obtain the worst-case point in the inner maximization problem, another assumption about concavity should be made.

\begin{assumption} \label{asp3}
    Objective function $h(\theta,\gamma;\xi,\xi')$ is $\mu_h$-strongly concave in $\xi'$ .
\end{assumption}

If Assumption \ref{asp3} holds, we can know if $\xi^{\star}_{t}$ means the maximizer of function $h$ and $\hat{\xi}_t^*$ is a $\varepsilon$-maximizer in $t$ iteration and , which means $|h(\xi^{\star}_t)-h(\hat{\xi}^{\star}_t)|\leq\varepsilon$, then we can get 
\begin{equation}
	\|\xi^{\star}_t-\hat{\xi}^{\star}_t\|^2\leq \frac{2\varepsilon}{\mu_h},
\end{equation}
and in algorithm, we set $\epsilon$ with the value $2\varepsilon/{\mu_{h}}$

In the following three situations, Assumption \ref{asp3} holds.
    \begin{enumerate}
        \item [{\romannumeral1}] When $c(f(\theta;X');Y')$ is non-concave and Assumption \ref{asp2} holds, we have $\mu_h=(\gamma-L_{\xi\xi})$ when $\gamma>L_{\xi\xi}$;
        \item [{\romannumeral2}] When  $c(f(\theta;X');Y')$ is concave in $\xi'$, we have $\mu_h=\gamma$;
        \item [{\romannumeral3}] When $c(f(\theta;X');Y')$ is $\mu_c$-strongly concave in $\xi'$, we have $L_h=\mu_c+\gamma$. 
    \end{enumerate}
\begin{lemma} \label{lm1}
    If Assumptions \ref{asp1}-\ref{asp3} hold, and $h^{\star}(\theta,\gamma;\xi)$ is differentiable, let $L_{\xi\gamma}$ denotes the Lipschitz constant for $\nabla_{\xi'}h(\theta,\gamma;\xi,\xi')$ in $\gamma$, then we can get:
    \begin{equation}
        \|\xi^{\star}(\theta_1,\gamma_1)-\xi^{\star}(\theta_2,\gamma_2)\|\leq
        L_1\|\gamma_1-\gamma_2\|+L_2\|\theta_1-\theta_2\|,
    \end{equation}
    \begin{equation}
        \|\nabla_{\theta}h^{\star}(\theta_1,\gamma_1;\xi)-\nabla_{\theta}h^{\star}(\theta_2,\gamma_2;\xi)\|_{*}\leq
        L_3\|\gamma_1-\gamma_2\|+L_4\|\theta_1-\theta_2\|,
    \end{equation}
    \begin{equation}
        \|\nabla_{\gamma}h^{\star}(\theta_1,\gamma_1;\xi)-\nabla_{\gamma}h^{\star}(\theta_2,\gamma_2;\xi)\|_{*}\leq
        L_1\|\gamma_1-\gamma_2\|+L_2\|\theta_1-\theta_2\|,
    \end{equation}
    where $L_1=\frac{L_{\xi\gamma}}{\mu_h}$, $L_2=\frac{L_{\xi\theta}}{\mu_h}$, $L_3=\frac{L_{\theta\xi}L_{\xi\gamma}}{\mu_h}$ and $L_4=L_{\theta\theta}+\frac{L_{\theta\xi}L_{\xi\theta}}{\mu_h}$.
\end{lemma}

\begin{corollary} [Lipschitz Smoothness of Convex Multifunction] \label{col1}
    If the gradient of $h(\theta,\gamma;\xi,\xi')$ for $\theta$ and $\gamma$ satisfies the implicit function theorem, and let $\xi^{\star}(\theta,\gamma)$ denotes the implicit function, which is a multifunction. Then when $\xi^{\star}(\theta)$ is convex in $\theta$ and $\gamma$, there must exist $k_{\theta}$ and $k_{\gamma}$ such that
    \begin{equation}
         \|\xi^{\star}(\theta_1,\gamma_1)-\xi^{\star}(\theta_2,\gamma_2)\| =
        k_{\gamma}|\gamma_2-\gamma_1\|+k_{\theta}\|\theta_2-\theta_1\|
    \end{equation}
\end{corollary}

Finally, to prove the convergence of proposed algorithm, another assumption about $H(\theta,\gamma)$ should be made:

\begin{assumption} \label{asp4}
    Assume that $H(\theta,\gamma)$ and $h^{\star}(\theta,\gamma;\xi)$ satisfy
    \begin{equation}
        \mathbb{E}_Q[\|\nabla_{\theta}H(\theta,\gamma)-\nabla_{\theta}h^{\star}(\theta,\gamma;\xi)\|^2]\leq \sigma_{\theta}^2,
    \end{equation}
    \begin{equation}
        \mathbb{E}_Q[\|\nabla_{\gamma}H(\theta,\gamma)-\nabla_{\gamma}h^{\star}(\theta,\gamma;\xi)\|^2]\leq \sigma_{\gamma}^2,
    \end{equation}
    \begin{equation}
        H(\theta_0,\gamma_0)-\inf_{\theta,\gamma}H(\theta,\gamma)\leq \Lambda,
    \end{equation}
    where $\Lambda>0$ is a constant
\end{assumption}

\begin{theorem}\label{thm1} (Algorithm Convergence)
    Let Assumptions \ref{asp1}-\ref{asp4} hold and take constant stepsize $\alpha=\beta=\sqrt{\frac{\Lambda}{L(\sigma_{\theta}^2+\sigma_{\gamma}^2)T}}$ and $\varepsilon=\frac{\mu_h}{(L_{\theta\xi}^2+1)\sqrt{T}}$, where $L^2=\max(L_1^2+L_3^2+L_1L_2+L_3L_4,L_2^2+L_4^2+L_1L_2+L_3L_4)$, Algorithm \ref{algorithm1} satisfies
    \begin{equation}
        \frac{1}{T}\sum_{t=0}^{T-1}\mathbb{E}_Q[\|\nabla_{\theta}H(\theta_t,\gamma_t)\|^2_2+\|\nabla_{\gamma}H(\theta_t,\gamma_t)\|_2^2]\leq (4\sqrt{L(\sigma_{\theta}^2+\sigma_{\gamma}^2)\Lambda} +2)\sqrt{\frac{1}{T}} +
        2\sqrt{\frac{L\Lambda}{\sigma_{\theta}^2+\sigma_{\gamma}^2}}\frac{1}{T}
    \end{equation}
\end{theorem}

Algorithm \ref{algorithm1} provides a $O(1/\sqrt{T})$ bound for expected gradient. When iteration $T$ tends to $\infty$, expected average gradient for $\theta$ and $\gamma$ both approach to zero with convergence speed no lower than $\sqrt{T}$, which can guarantee Algorithm \ref{algorithm1} find a stationary point in non-convex optimization problems and optimal point in convex optimization problems. When objective function $c(f(\theta;X);Y)$ is convex in $\theta$, Corollary \ref{col2} provides a formal bound on the distance between $\theta_T$ and $\theta_{R}$.

\begin{corollary}[Solution Convergence] \label{col2}
    If $c(f(\theta;X);Y)$ is $\lambda$-strongly convex in $\theta$ and assume $\mathbb{E}_{Q}[\|\theta_0-\theta_R\|] \leq M_0$, then a formal bound can be provided as follows:
    \begin{equation}
        \mathbb{E}_{Q}[\|\theta_T-\theta_R\|^2]\leq O(\frac{1}{\sqrt{T}})
    \end{equation}
    where $\theta_T$ is the algorithm solution after $T$ iterations, and $\theta_R$ is the robust optimal solution.
\end{corollary}

\begin{corollary}[Objective Value Gap]\label{col3}
    If Assumption \ref{asp2} holds, based on Corollary \ref{col2}, the gap about objective value between algorithm result and optimal result can be bounded as:
    \begin{equation}
        0 \leq H(\theta_t,\gamma_t) - H(\theta^{\star},\gamma^{\star}) \leq \frac{L_{\theta\theta}}{2}\|\theta_t-\theta^{\star}\| \leq O(\frac{1}{\sqrt{T}}),
    \end{equation}
    which means it also has a $O(\frac{1}{\sqrt{T}})$ bound.
\end{corollary}

Corollaries \ref{col2} and \ref{col3} show that the solutions and objective values of the algorithm are closer to the optimal ones to the DRO model \eqref{e4}. It also provides a guideline for the practitioners to trade off between computational efficiency and solution accuracy. For different levels of tolerance, we can use the corollaries to determine the number of iterations to run the algorithm.

With $r_{algorithm}$, the rest of performance gap of the proposed algorithm is determined by $r_{DRO}$. Without losing generality, if a DRO model \eqref{e5} is designed to have $r_{DRO}=\nu$, then the total gap is $O(\frac{1}{\sqrt{T}})+\nu$. In the following, we exemplify the gap in two special cases. The first is the one where the nominal distribution is the empirical distribution based on the historical samples, in which $r_{DRO}$ is bounded by the size of dataset $N$. The second one is the online case where the sample can be repeatedly drawn from the oracle distribution, in which $r_{DRO}$ is bounded by training iterations $T$.

\subsection{Special Case 1: Empirical Nominal Distribution}\label{sec4.2}
In this section, we consider an empirical nominal distribution, with dataset $\{\xi_i\}_{i=1}^N$.
Many literature has shown that when the radius $\rho$ is properly set according to the number of data $N$, one can guarantee the out-of-sample performance of the DRO model with empirical nominal distributions.
For example, under different assumptions and for a given confidential level that the true distribution is included in the ambiguity set, denoted as $q$, \cite{mohajerin2018data, fournier2015rate} provide a lower bound on the minimal radius required, which is of the order $O(N^{-1/s})$, where $s$ is the dimension of the random variable, and  provide $O(N^{-1/2})$ order is proved by \cite{gao2023finite} and \cite{ji2021data}. To make the paper more self-included, we present some similar results with more general assumptions and smaller minimal required radius.

\begin{assumption} 
    Assume there exist some $\xi_0\in\Xi$ such that the exponential integrability condition (defined by \cite{wang2010sanov})
	\begin{equation}
		\Lambda(a):=\log{\int_{\Xi} \exp\{a d(\xi,\xi_0)\}\,dQ(\xi)}<+\infty,\quad \forall a>0 \notag
	\end{equation}
	holds true for all $a>0$.
\end{assumption}

\begin{theorem}[Probability Guarantee] \label{probguarantee}
    Assume assumption \ref{asp1} and \ref{asp2} hold, and the Wasserstein ball radius $\rho_N$ is smaller than support space diameter $D$, then A lower bound on the probability that the Wasserstein distance between the empirical distribution $P_N$ and the real distribution $P_{true}$ does not exceed $\rho_N$ is given by
	\begin{equation} \label{e6}
		P(W_d(P_{true},P_N)\leq\rho_N)\geq
		\left\{
		\begin{aligned}
			&1-e^{-N\frac{\rho_N^2}{2(1+D^2)}}\quad &&if\; &&D\geq 1,\\
			&1-e^{-N\frac{\rho_N^2}{4D^2}}&&if &&D < 1.\\
		\end{aligned}
		\right.
	\end{equation}
\end{theorem}

In Theorem \ref{probguarantee}, the probability of ambiguity containing real distribution can be calculated, once the radius $\rho$ is determined. With the same radius $\rho_N$, our result provides a higher probability than that in \cite{ji2021data}, which means a smaller radius on the same confidence level. Compared with the result in \cite{gao2023finite}, less assumption is needed here and the radius is easier to be calculated in Theorem \ref{probguarantee}. Based on the result in the above theorem, if the confidence level $q$ is determined, the radius can also be deduced.

\begin{corollary}[Radius with High Probability] \label{radiuschosen}
    Assume assumption \ref{asp1}, \ref{asp2} and $\rho<D$ hold, then if
	\begin{equation}
		\rho_N\geq
		\left\{
		\begin{aligned} \label{e12}
			&\sqrt{-\frac{1}{N}\log (1-q)[2(1+D^2)]}\quad &&if\; &&D\geq 1,\\
			&2D\sqrt{-\frac{1}{N}\log(1-q)}&&if &&D < 1.\\
		\end{aligned}
		\right.
	\end{equation}
	then
	\begin{align}
		P(W_d(P_{true},P_N)\leq\rho_N)\geq q \quad \Longrightarrow \quad  P_{true}\in\mathcal{P}\  with\  probability\ at\ least\ q. \notag
	\end{align}
\end{corollary}

Corollary \ref{radiuschosen} provides a method to choose the radius $\rho$ with high confidence level $q$. After determining the radius of Wasserstein ball, gap $r_{DRO}$ can be bounded by data number $N$, which is also known as finite sample guarantee. Let $C_N(\theta)=\mathbb{E}_{P_N}[c(f(\theta;X);Y)]$ denotes the objective with empirical distribution, then $r_{DRO}$ can be calculated by:
\begin{align}
    r_{DRO} &= C(\theta_R)-C(\theta^{\star}) \notag\\
        &= [C(\theta_R)-C_{DRO}(\theta_R;\rho_N)]+[C_{DRO}(\theta_R;\rho_N)-C_{DRO}(\theta^{\star};\rho_N)] \notag\\
        &\ \ +[C_{DRO}(\theta^{\star};\rho_N)-C_N(\theta^{\star})]+[C_N(\theta^{\star})-C(\theta^{\star})].
\end{align}

\begin{theorem}[DRO Gap with Empirical Distribution] \label{finiesampleguarantee}
    If Assumption \ref{asp1} holds, then With probability at least $q(1-e^{-t})$,
    \begin{equation}
        r_{DRO} \leq \sqrt{\frac{\tau t}{N}}L_c + 0 + \rho_N L_c + O(\sqrt{\frac{1}{N}}),
    \end{equation}
    where $\tau=2(1+D^2)$ when $D\geq1$ and $\tau=4$ when $D<1$, $c(f(\theta;X);Y)$ is $L_c$-Lipschitz continue in $\xi$, and $\rho_N$ is the radius decided in Theorem \ref{probguarantee} with probability $q$. Thus, the high probability upper bound of $r_{DRO}$ is $O(1/\sqrt{N})$. 
\end{theorem}

\begin{corollary}[Finite Sample Guarantee] \label{finitesampleerror}
    The gap $r$ is bounded by $O(1/\sqrt{N}+1/\sqrt{T})$ when nominal distribution $Q$ is empirical distribution $P_N$.
\end{corollary}

\begin{corollary}[Out-of-sample Guarantee] \label{outofsamplebound}
    Consider the generalization bound (defined in \cite{shafieezadeh2019regularization} Theorem 36), which is also known as out-of-sample-bound with the form as:
    \begin{equation}
    |C_{DRO}(\theta_T;\rho_N)-C(\theta_T)| \leq \rho_N L_c + O(\frac{1}{\sqrt{N}}) \leq O(\frac{1}{\sqrt{N}}).
    \end{equation}
\end{corollary}

\subsection{Special Case 2: Online Learning Problems} \label{sec4.3}
As we have discussed, our algorithm free the construction of the DRO model with specific knowledge on the nominal distribution as long as we can sample from it. Online learning is one of the specific instances of such a case, where the nominal distribution is the ground-truth distribution, and new samples from the distribution can be taken as the bootstrapped ones in our algorithm. 

However, different from the DRO model with empirical distribution, the dataset in online learning is updated iteration by iteration which results in an increasing data number $N$. Thus, a fixed radius $\rho$ should be chosen without the knowledge of data size. From Theorem \ref{probguarantee}, we can know, the probability that ambiguity set contains ground-truth distribution will increase with the $N$ increasing.

After deciding the radius, We can use the results of the algorithm gap to determine the online regret when applying our algorithm in the setting. For online learning, its regret is defined as
\begin{equation}
    Regret_T(\theta^{\star},\gamma^{\star}) = \sum_{t=1}^T \big(h(\theta_t,\gamma_t;\xi_t,\xi_t^{\star})-h(\theta^{\star},\gamma^{\star};\xi_t,\xi^{\star}) \big).
\end{equation}

If Algorithm \ref{algorithm1} can only generate a stationary solution, the average regret of online learning will never converge to zero with iterations tends to $\infty$. Therefore, it is a reasonable assumption that $h(\theta,\gamma;\xi,\xi')$ is convex in $\theta$, when analyzing online learning regret.

\begin{corollary} [Online Regret]\label{regret}
    If $h(\theta,\gamma;\xi,\xi')$ is convex in $\theta$, then the expectation of regret satisfies
    \begin{equation}
        \mathbb{E}_Q [Regret_T(\theta^{\star},\gamma^{\star})] \leq O(\sqrt{T}),
    \end{equation}
    and it is a sub-linear function about $T$. The speed of convergence is $O(\frac{1}{\sqrt{T}})$.
\end{corollary}

In real life, the ground-truth distribution of random variables may change over time. From the Corollary \ref{regret}, the convergence speed of Algorithm \ref{algorithm1} is $1/\sqrt{T}$. This implies that when the distribution changes, new proposed algorithm can capture this change and obtain the optimal solution of latest distribution with the speed of $1/\sqrt{T}$.

\section{Numerical Experiments}\label{sec5}
We validate both effectiveness and computational efficiency of our proposed method through data-driven newsvendor problem, a common benchmark of many contextual optimization literature, e.g., \citep{ban2019big,zhang2023optimal}. We compare our proposed DA-SGD algorithm against several benchmarks, including empirical risk minimization decision policy without and with $l^1$ regularization (NV-ERM1 and NV-ERM2) in \cite{ban2019big}, conditional stochastic optimization using random forests (RF) in \cite{bertsimas2020predictive}, robust policy and its Shapley policy (Shapley) in \cite{zhang2023optimal}.

All experiments are performed in Windows 11 using Python 3.10 with an optimization solver Gurobi 9.5.2 on a Dell Inspiron 15 7510 with Intel(R) Core(TM) i5-11400H CPU and 16 GB RAM.

\subsection{Experiment Settings} \label{sec5.1}
Newsvendor model is a classical problem in inventory control, which has the objective function as follows:
\begin{equation}
    c(z;Y)=C_b(Y-z)^++C_h(z-Y)^+
\end{equation}
\noindent 
where $(\cdot)^+$ means the non-negative part, $Y$ and $z$ denotes the demand and order of the product, and $C_b$ and $C_h$ are, respectively, the unit back-ordering and holding costs.

Thus, the contextual optimization problem of Newsvendor should have the following form:
\begin{equation} \label{e21}
    c(f(\theta;X);Y)=C_b(Y-f(\theta;X))^++C_h(f(\theta;X)-Y)^+
\end{equation}

Unfortunately, \eqref{e21} does not satisfy Assumption \ref{asp2}. In order to solve the Newsvendor model with DA-SGD algorithm, a quadratic function is used to approximate the inflection point in the small neighbor $(Y-\delta,Y+\delta)$, which is:
\begin{align} \label{e22}
    c_{new}(f(\theta;X);Y)=
    \left\{
\begin{aligned}
&C_b(Y-f(\theta;X)) \quad &&if \ f(\theta;X) \leq Y-\delta \\ 
&a_1 (f(\theta;X))^2 +a_2 f(\theta;X) + a_3   &&if Y-\delta<f(\theta;X)<Y+\delta\\
&C_h(f(\theta;X)-Y) \quad &&if \ f(\theta;X) \geq Y+\delta\\
\end{aligned}
\right.
\end{align}
where $a_1=(C_b+C_h)/4\delta$, $a_2=-(C_b(Y+\delta)+C_h(Y-\delta))/2\delta$, $a_3=C_b\delta+(C_b(Y+3\delta)+C_h(Y-\delta))(Y-\delta)/4\delta$, and $\delta$ is the radius of the approximation neighbor, which has the effect on optimization error and limitation of Lagrangian factor $\lambda$. The reason and benefits of choosing such an approximation function is:
\begin{itemize}
    \item $c_{new}(f(\theta;X);Y-\delta)=c(f(\theta;X);Y-\delta)$;
    \item $c_{new}(f(\theta;X);Y+\delta)=c(f(\theta;X);Y+\delta)$;
    \item $\nabla_{\xi} c_{new}(f(\theta;X);Y-\delta)=\nabla_{\xi} c(f(\theta;X);Y-\delta)$;
    \item $\nabla_{\xi} c_{new}(f(\theta;X);Y+\delta)=\nabla_{\xi} c(f(\theta;X);Y+\delta)$.
\end{itemize}

These four properties guarantee the continuity and smoothness of new objective function, which makes Assumption \ref{asp2} hold. Additionally, based on these properties, if we consider a linear decision policy $f(\theta;X)=\theta^{\top}x$, then \eqref{e22} becomes a smooth function with its gradient function satisfying $(C_b+C_h)\theta^2/2 \delta$-Lipschitz. The largest decision error between \eqref{e22} and \eqref{e21} is less than $(C_b+C_h)\delta/4$. The size of approximation neighbour can be chosen to balance Lipschitz constant and largest error, and in experiment it is chosen based on several pre-experiment.

We adapt the experiment setup from \cite{zhu2012semiparametric} and consider a linear demand model with additive noise:
\begin{equation} \label{linearfunction}
    y=\theta^\mathsf{T}x+\theta_0+\epsilon,
\end{equation}
where $\theta$ is an $s$-dimensional vector representing the ground-truth coefficients, with $s\in \{100,1000\}$. The feature variables $x\in \mathcal{R}^{s}$  are drawn from an uniform distribution from $[0,1]^s$ (or from multivariate Gaussian distribution with mean and covariance). The noise term follows a Gaussian distribution $\epsilon \sim \mathcal{N}(0,\sigma^2_{\epsilon})$, with $\sigma_{\epsilon}^2\in \{0.1, 0.5,1,5\}$.

To enlarge the effect of unbalanced prediction error, we set $C_h=0.2$ and $C_b=1$. To study the impact of the sample size $n$, especially for the small data situations, we vary the number of data points $n\in\{10, 30, 50, 70, 100,300, 500\}$. To avoid the coincidence, the testing data size is fixed as $10,000$, which is much larger than training data size. We run 20 repeated experiments.

\subsection{Experiment Result and Analysis}
\begin{figure}[htbp]
    \centering
    {\includegraphics[scale=0.5]{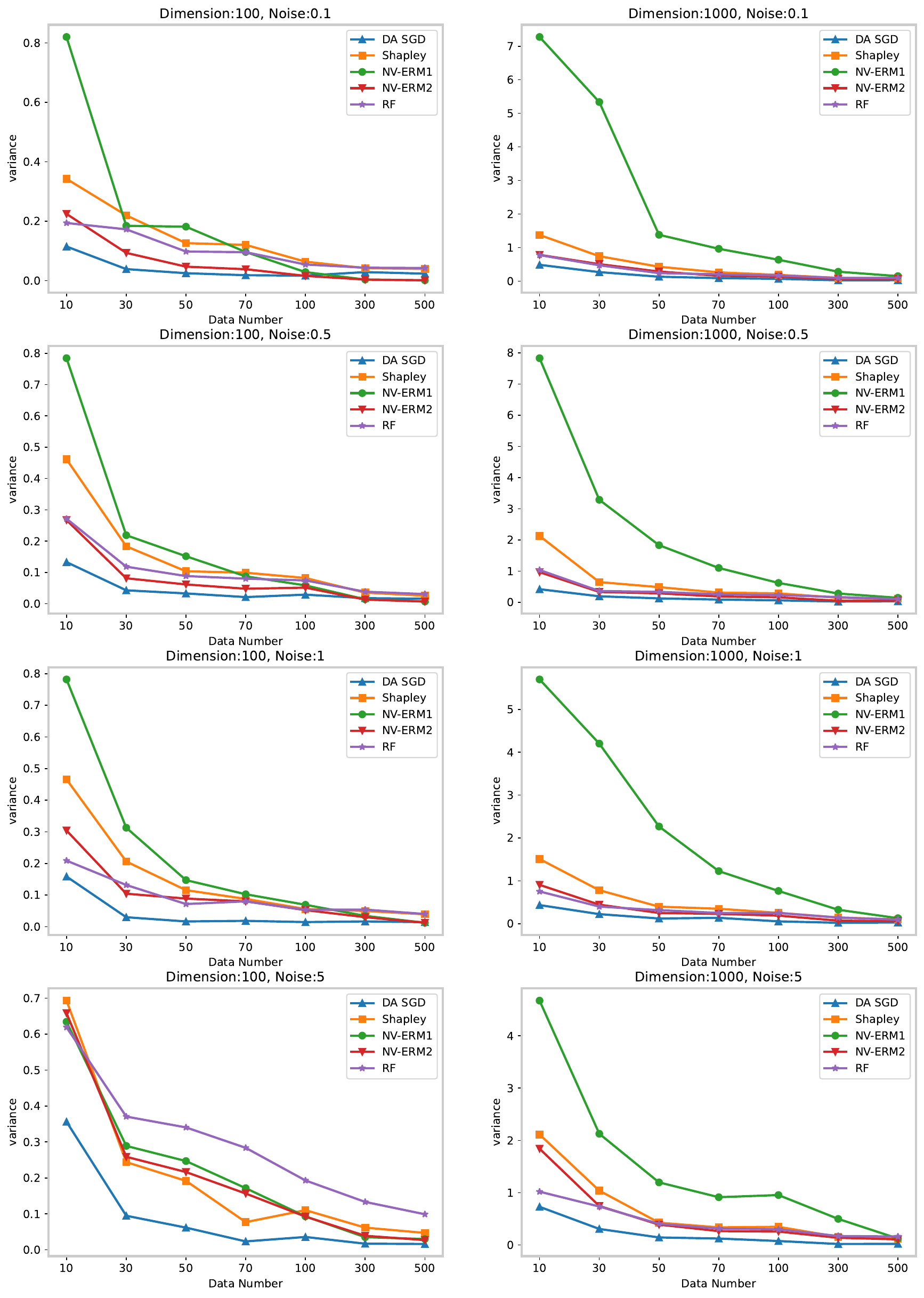}}
    \caption{Variance of out-of-sample objective over 20 instances}
    \label{linearvariance}
\end{figure}

In this section, we apply linear function \eqref{linearfunction} to generate the demand, with noise range in $\sigma\in\{0.1, 0.5,1,5\}$ standing for the quality of dataset, run 20 repeated experiments and draw the line chart in Figure \ref{linearvalue} and \ref{linearvariance}. In each experiment, the out-of-sample objective values of our method and benchmarks are calculated with the same training and testing data set. The line in Figure \ref{linearvalue} represents the average expected objective value over 20 rounds, while the line in Figure \ref{linearvariance} indicates the variance of each methods. We have the following observations:

\begin{figure}[htbp]
    \centering
    {\includegraphics[scale=0.5]{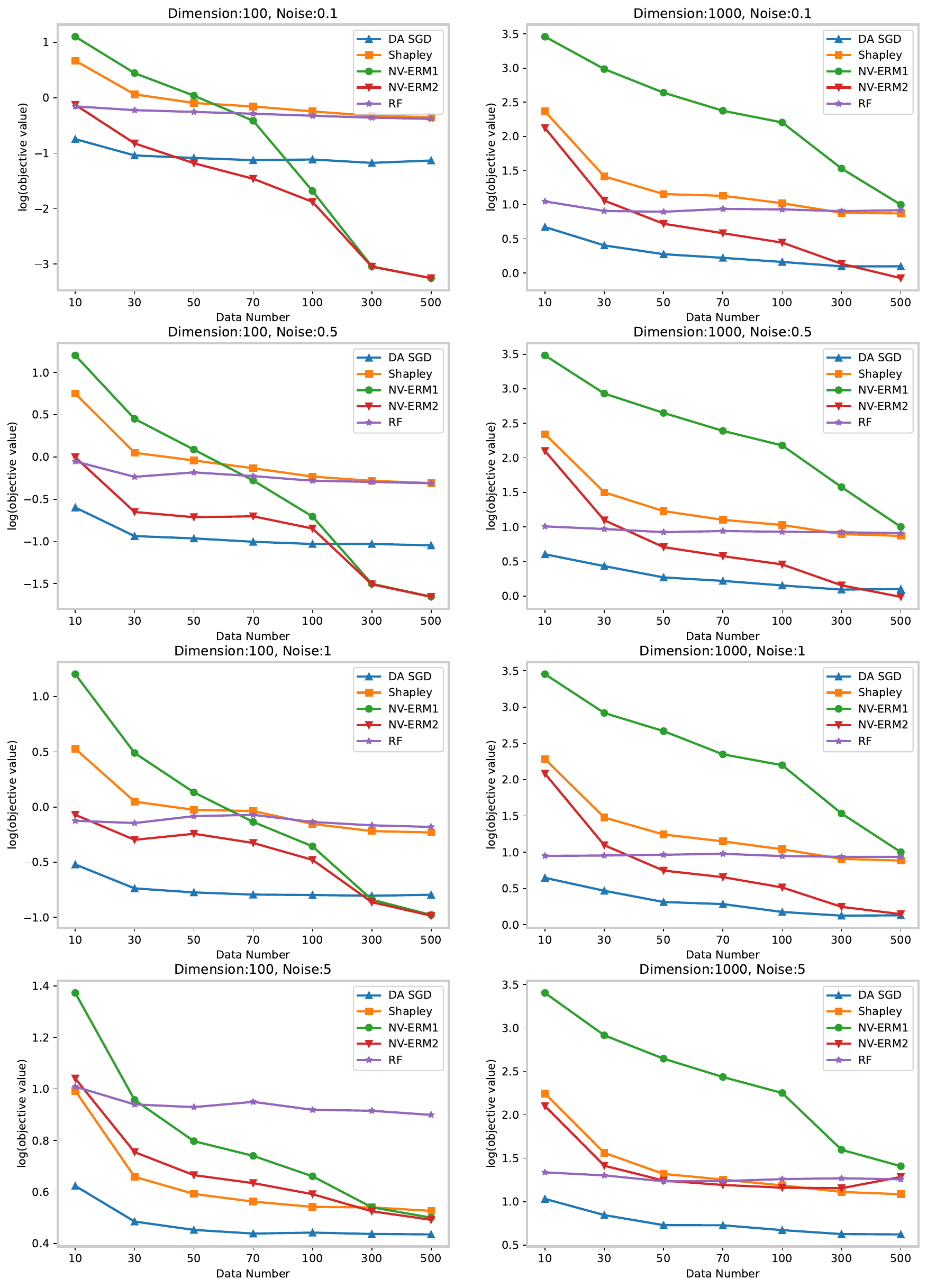}}
    \caption{Average out-of-sample objective over 20 instances}
    \label{linearvalue}
\end{figure}

Figure \ref{linearvalue} illustrates the log-transformed average objective results over 20 runs. The Figure shows that when data quality is high (small $\sigma$), DA-SGD outperforms other benchmarks with small sample size. As data size increases, NV-ERM methods surpass DA-SGD. Shapley and RF exhibit poor performance under current experiment settings. When data quality deteriorates, the advantage of DA-SGD becomes more obvious even with large data set. When dimension comes to 1000, the small-data-large-scale phenomenon becomes more evident, and DA-SGD outperforms other benchmarks. Thus, from this result, DA-SGD is least affected by the size of data set, and exhibit a relatively smooth trend as size increasing. On the contrary, NV-ERM methods are highly influenced by data volume, performing worst with small data size. DA-SGD method demonstrates stronger robustness, performing better with small, low-quality data.

Figure \ref{linearvariance} shows the variance across 20 experiment runs. The overall fluctuation of the DA-SGD method is smallest among these methods, while NV-ERM methods have the largest variance with small data and is significantly influenced by the randomness of data generation. As data quality deteriorates, the variances of all methods increase, but the impact of DA-SGD is smaller than other benchmarks. This indicates that the DA-SGD method is highly stable, minimally affected by data fluctuations, and exhibits strong robustness.

DA-SGD and Shapley solve the DRO problems, while the other benchmarks solve initial problems directly, which is much easier to solve than DRO problems. Thus, We only report the average computational time of DA-SGD and Shapley in Table \ref{lineartime}, in which the computational time contains both training time and testing time. Because noise $\sigma$ has little effect on the computational time, the average computational time is calculated over all 20 experiments and all different $\sigma$. Based on the result, we have the following observations:

\begin{table}
    \centering
    \begin{tabular}{ccrr}
    \toprule
    Dimension ($s$) & Data Number ($d$) & DA-SGD & Shapley\\
    \toprule
    100  & 10 & 1.5558 & 14.6389\\
         & 30 & 1.5158 & 37.9416\\
         & 50 & 1.5066 & 62.1709\\
         & 70 & 1.5048 & 86.1249\\
         & 100& 1.5113 & 125.1831\\
         & 300& 1.4911 & 434.9086\\
         & 500& 1.4865 & 574.5781\\
    1000 & 10 & 1.8482 & 73.4890 \\
         & 30 & 1.8425 & 214.5982\\
         & 50 & 1.8534 & 357.0431\\
         & 70 & 1.8562 & 499.2080\\
         & 100& 1.8418 & 712.4738\\
         & 300& 1.8375 &2188.1067\\
         & 500& 1.8376 &3494.8938\\
    \bottomrule
    \end{tabular}
    \caption{Average Computational Time (in Seconds) per DRO Problem Instance}
    \label{lineartime}
\end{table}

As mentioned before, DA-SGD is $O(T\cdot K)$ and its computational time is mainly influenced by training parameters $T$ and $K$, while the size of data set can only have impact on solving time through increasing the difficulty of matrix calculation. Thus, with different size of data set, the time of DA-SGD is stable. However, the computational time of Shapley is $O(n)$ in training and $O(n^2)$ in extension, which makes it much longer then DA-SGD and more sensitive to the data size. 

\section{Conclusions}\label{sec6}
This paper proposes an efficient approach to finding a robust, optimal, end-to-end policy for the feature-based stochastic optimization problems. This new proposed algorithm utilizes data augmentation method to solve the inner maximization problems of the duality problem of DRO, and then applies the worst-case samples to train the parameters by SGD method, which can be extended into online learning. In this paper, some performance guarantees are provided for DA-SGD algorithm. For arbitrary nominal distribution, there will be a $O(1/\sqrt{T})+\nu$ bound for the gap between optimal solution and algorithm solution, while in two special cases, it will become $O(1/\sqrt{T}+1/\sqrt{N})$ with nominal empirical distribution and $O(1/\sqrt{T})$ in online learning problems.

There are three potential future research. Firstly, because most stochastic optimization methods are suitable for unconstrained problems, it would be worthwhile to extend it to constrained stochastic optimization problems. Secondly, this algorithm is used to solve DRO problems to overcome small data challenges, and it is interesting to consider more data challenges, such as data missing. Thirdly, it is valuable if more additional factors that influence decision-making can be combined with optimization, for example using bi-level optimization to select feature.

\bibliographystyle{informs2014}
\bibliography{DA-OPT}

\clearpage
\begin{APPENDICES}
    \section{Proof of Section \ref{sec4.1}} \label{secA}
    \subsection{Proof of Lemma \ref{lm1}} \label{secA.1}
    When Assumption \ref{asp3} holds, $h(\theta,\gamma;\xi,\xi')$ is $\mu_h$-strongly concave in $\xi'$, and $\xi^{\star}(\theta,\gamma)$ is the worst-case point with $\theta$ and $\gamma$. For any $\theta$ and $\gamma$, optimality of $\xi^{\star}(\theta,\gamma)$ implies that $\nabla_{\xi'}h(\theta,\gamma;\xi,\xi')^{\top}(\xi'-\xi^{\star}(\theta,\gamma))\leq0$. Then, based on strongly concavity, let $\xi^{\star}_1=\xi^{\star}(\theta_1,\gamma_1)$ and $\xi^{\star}_2=\xi^{\star}(\theta_2,\gamma_2)$, we have
    \begin{equation}
        h(\theta_2,\gamma_2;\xi,\xi^{\star}_2) \leq h(\theta_2,\gamma_2;\xi,\xi^{\star}_1) + \nabla_{\xi'}h(\theta_2,\gamma_2;\xi,\xi^{\star}_1)^{\top}(\xi^{\star}_2-\xi^{\star}_1)-\frac{\mu_h}{2}\|\xi^{\star}_2-\xi^{\star}_1\|^{2}, \notag
    \end{equation}
    \begin{align}
        h(\theta_2,\gamma_2;\xi,\xi^{\star}_1) &\leq h(\theta_2,\gamma_2;\xi,\xi^{\star}_2) + \nabla_{\xi'}h(\theta_2,\gamma_2;\xi,\xi^{\star}_2)^{\top}(\xi^{\star}_1-\xi^{\star}_2)-\frac{\mu_h}{2}\|\xi^{\star}_2-\xi^{\star}_1\|^{2} \notag \\
        &\leq h(\theta_2,\gamma_2;\xi,\xi^{\star}_2) -\frac{\mu_h}{2}\|\xi^{\star}_2-\xi^{\star}_1\|^{2}. \notag
    \end{align}    
    Summing these two inequalities gives
    \begin{align}
        \mu_h \|\xi^{\star}_2-\xi^{\star}_1\|^{2} &\leq  \nabla_{\xi'}h(\theta_2,\gamma_2;\xi,\xi^{\star}_1)^{\top}(\xi^{\star}_2-\xi^{\star}_1) \notag \\
        &\leq  [\nabla_{\xi'}h(\theta_2,\gamma_2;\xi,\xi^{\star}_1)-
        \nabla_{\xi'}h(\theta_1,\gamma_2;\xi,\xi^{\star}_1)+
        \nabla_{\xi'}h(\theta_1,\gamma_2;\xi,\xi^{\star}_1)-
        \nabla_{\xi'}h(\theta_1,\gamma_1;\xi,\xi^{\star}_1)
        ]^{\top}(\xi^{\star}_2-\xi^{\star}_1), \notag
    \end{align}
    where $\nabla_{\xi'}h(\theta_1,\gamma_1;\xi,\xi^{\star}_1)
        ^{\top}(\xi^{\star}_2-\xi^{\star}_1)\leq0$. Using Holder's inequality, we can get
    \begin{align}
        \mu_h \|\xi^{\star}_2-\xi^{\star}_1\|^{2} &\leq
        [\|\nabla_{\xi'}h(\theta_2,\gamma_2;\xi,\xi^{\star}_1)-\nabla_{\xi'}h(\theta_1,\gamma_2;\xi,\xi^{\star}_1)\|_{\star} \notag \\ 
        &\ \ +\|\nabla_{\xi'}h(\theta_1,\gamma_2;\xi, \xi^{\star}_1) -\nabla_{\xi'} h(\theta_1,\gamma_1;\xi, \xi^{\star}_1)\|_{\star}]\|\xi^{\star}_2-\xi^{\star}_1\|, \notag
    \end{align}
    which is,
    \begin{equation}
        \|\xi^{\star}_2-\xi^{\star}_1\|^{2} \leq \frac{L_{\xi\gamma}\|\gamma_1-\gamma_2\| + L_{\xi\theta}\|\theta_2-\theta_1\|}{\mu_h}. \notag
    \end{equation}
    To see the other two inequalities, we should know that 
    \begin{equation}
        \nabla_{\theta}h(\theta,\gamma;\xi,\xi') = \nabla_{\theta}c(f(\theta;x');y'), \notag    
    \end{equation}
     \begin{equation}
        \nabla_{\gamma}h(\theta,\gamma;\xi,\xi') = \rho - d(\xi,\xi'). \notag     
     \end{equation}
    Then to calculate the gap between gradient with worst-case point can be,
    \begin{align}
        \|\nabla_{\theta}h(\theta_2,\gamma_2; \xi,\xi^{\star}_2) - \nabla_{\theta}h(\theta_1, \gamma_1; \xi,\xi^{\star}_1)\|_{\star} &\leq \|\nabla_{\theta}h(\theta_2,\gamma_2; \xi,\xi^{\star}_2) - \nabla_{\theta}h(\theta_2, \gamma_2; \xi,\xi^{\star}_1)\|_{\star}  \notag\\ &\ \ +\|\nabla_{\theta}h(\theta_2,\gamma_2; \xi,\xi^{\star}_1) - \nabla_{\theta}h(\theta_1,\gamma_1; \xi,\xi^{\star}_1)\|_{\star} \notag,
    \end{align}
    for the first part, we have
    \begin{equation}
        \|\nabla_{\theta}h(\theta_2,\gamma_2; \xi,\xi^{\star}_2) - \nabla_{\theta}h(\theta_2, \gamma_2; \xi,\xi^{\star}_1)\|_{\star} \leq L_{\theta\xi}\|\xi^{\star}_2-\xi^{\star}_1\|, \notag
    \end{equation}
    and for the second one, it should be
    \begin{align}
        \|\nabla_{\theta}h(\theta_2,\gamma_2; \xi,\xi^{\star}_1) - \nabla_{\theta}h(\theta_1,\gamma_1; \xi,\xi^{\star}_1)\|_{\star} &\leq 
        \|\nabla_{\theta}h(\theta_2,\gamma_2; \xi,\xi^{\star}_1) - \nabla_{\theta}h(\theta_1,\gamma_2; \xi,\xi^{\star}_1)\|_{\star} \notag \\ &\ \ + \|\nabla_{\theta}h(\theta_1,\gamma_2; \xi,\xi^{\star}_1) - \nabla_{\theta}h(\theta_1,\gamma_1; \xi,\xi^{\star}_1)\|_{\star} \notag \\
        &\leq L_{\theta\theta}\|\theta_2-\theta_1\|. \notag
    \end{align}
    Thus, 
    \begin{equation}
        \|\nabla_{\theta}h(\theta_2,\gamma_2; \xi,\xi^{\star}_2) - \nabla_{\theta}h(\theta_1, \gamma_1; \xi,\xi^{\star}_1)\|_{\star} \leq L_{\theta\xi}\|\xi^{\star}_2-\xi^{\star}_1\| + L_{\theta\theta}\|\theta_2-\theta_1\|. \notag
    \end{equation}
    It is the same to calculate the gap of gradient about $\gamma$,
    \begin{align}
        \|\nabla_{\gamma}h(\theta_2,\gamma_2; \xi,\xi^{\star}_2) - \nabla_{\gamma}h(\theta_1, \gamma_1; \xi,\xi^{\star}_1)\|_{\star} &\leq \|d(\xi,\xi^{\star}_2)-d(\xi,\xi^{\star}_1)\|_{\star} \notag \\ 
        &\leq \|d(\xi^{\star}_2,\xi^{\star}_1)\|_{\star} \notag \\
        &\leq \|\xi^{\star}_2-\xi^{\star}_1\|. \notag
    \end{align}
    These are the desired results in Lemma \ref{lm1}.

    \begin{lemma} \label{smoothness}
        If Assumption \ref{asp3} holds, then based on Lemma \ref{lm1}, we have
        \begin{equation}
            \|\nabla H(\theta_2,\gamma_2)-\nabla H(\theta_1,\gamma_1)\|^2 \leq L^2\|(\theta_2,\gamma_2)-(\theta_1,\gamma_1)\|^2, \notag
        \end{equation}
    where $\nabla H$ denotes the gradient of $H(\theta,\gamma)$ in both $\theta$ and $\gamma$, and $L^2=\max(L_1^2+L_3^2+L_1L_2+L3L_4,L_2^2+L_4^2+L_1L_2+L_3L_4)$. This result shows the $L-$Lipschitz smoothness of $H(\theta,\gamma)$.
    \end{lemma}

    Lemma \ref{smoothness} is easy to prove, with the result obtained in Lemma \ref{lm1},
    \begin{align}
        \|\nabla H(\theta_2,\gamma_2)-\nabla H(\theta_1,\gamma_1)\|^2 &= \|\nabla_{\theta} H(\theta_2,\gamma_2)-\nabla_{\theta} H(\theta_1,\gamma_1)\|^2+\|\nabla_{\gamma} H(\theta_2,\gamma_2)-\nabla_{\gamma} H(\theta_1,\gamma_1)\|^2 \notag \\
        &\leq (L_1^2+L_3^2)\|\gamma_2-\gamma_1\|^2 + 2(L_1L_2+L_3L_4)\|\gamma_2-\gamma_1\|\|\theta_2-\theta_1\|+(L_2^2+L_4^2)\|\theta_2-\theta_1\|^2 \notag \\
        &\leq  (L_1^2+L_3^2+L_1L_2+L_3L_4)\|\gamma_2-\gamma_1\|^2 + (L_2^2+L_4^2+L_1L_2+L_3L_4)\|\theta_2-\theta_1\|^2 \notag \\
        &\leq L^2\|(\theta_2,\gamma_2)-(\theta_1,\gamma_1)\|^2. \notag
    \end{align}
    
    \subsection{Proof of Corollary \ref{col1}} \label{secA.2}
    If implicit multifunction $\xi^{\star}(\theta,\gamma)$ is convex in $\theta$ and $\gamma$, then
    \begin{equation}
        \xi^{\star}(\alpha\theta_1+(1-\alpha)\theta_2,\gamma) \supset \alpha\xi^{\star}(\theta_1,\gamma) + (1-\alpha) \xi^{\star}(\theta_2,\gamma) \quad
        \forall \theta_1, \theta_2 \ and \ \alpha\in[0,1] \notag
    \end{equation}
    \begin{equation}
        \xi^{\star}(\theta,\beta\gamma_1+(1-\beta)\gamma_2) \supset \beta\xi^{\star}(\theta,\gamma_1) + (1-\beta) \xi^{\star}(\theta,\gamma_2) \quad
        \forall \gamma_1, \gamma_2 \ and \ \beta\in[0,1] \notag
    \end{equation}
    And for each fixed $\theta$ and $\gamma$, only one $\xi^{\star}$ can get. Let $\theta_3=\alpha\theta_1+(1-\alpha)\theta_2$ and $\theta_2>\theta_1$, then we have 
    \begin{align}
        \alpha(\xi^{\star}(\theta_3,\gamma)-\xi^{\star}(\theta_1,\gamma)) = (1-\alpha)(\xi^{\star}(\theta_2,\gamma)-\xi^{\star}(\theta_3,\gamma)). \notag
    \end{align}
    Both sides are divided by $\|\theta_2-\theta_1\|$,
    \begin{align}
        \frac{\alpha(\xi^{\star}(\theta_3,\gamma)-\xi^{\star}(\theta_1,\gamma))}{\|\theta_2-\theta_1\|} &= \frac{(1-\alpha)(\xi^{\star}(\theta_2,\gamma)-\xi^{\star}(\theta_3,\gamma))}{\|\theta_2-\theta_1\|} \notag\\
        \frac{\alpha(\xi^{\star}(\theta_3,\gamma)-\xi^{\star}(\theta_1,\gamma))}{\|\theta_3-\theta_1\|} &= \frac{(1-\alpha)(\xi^{\star}(\theta_2,\gamma)-\xi^{\star}(\theta_3,\gamma))}{\|\theta_2-\theta_3\|}=k_{\theta} \notag
    \end{align}
    The result above is same for $\gamma$, and let $k_{\gamma}$ denotes the ratio. Based on this, it is easy to get Corollary \ref{col1}.
    
    \subsection{Proof of Theorem \ref{thm1}} \label{secA.3}
    For convenience, let $g_{\theta}^t=\nabla_{\theta}h(\theta_t,\gamma_t;\xi_t,\xi'_t)$,  $g_{\gamma}^t=\nabla_{\gamma}h(\theta_t,\gamma_t;\xi_t,\xi'_t)$ and $\hat{\xi}^{\star}_t$ denotes the $\varepsilon$-maximizer of $h(\theta_t,\gamma_t;\xi_t,\xi'_t)$, then based on SGD method, we have
    \begin{equation}
        \theta_{t+1} = \theta_{t} - \alpha_tg_{\theta}^t, \quad\quad
        \gamma_{t+1} = \gamma_{t} - \beta_tg_{\gamma}^t. \notag
    \end{equation}
    By a Taylor expansion using the $L-$smoothness of the function $H(\theta,\gamma)$, which is shown in Lemma \ref{smoothness}, we have
    \begin{align}
        H(\theta_{t+1},\gamma_{t+1}) &\leq H(\theta_t,\gamma_t) + \langle \nabla_{\theta}H(\theta_t,\gamma_t), \theta_{t+1}-\theta_t \rangle + \frac{L}{2}\|\theta_{t+1}-\theta_{t}\|^2 + \langle \nabla_{\gamma}H(\theta_t,\gamma_t), \gamma_{t+1}-\gamma_t \rangle + \frac{L}{2}\|\gamma_{t+1}-\gamma_{t}\|^2 \notag \\ 
        &= H(\theta_t,\gamma_t) -\alpha_t\|\nabla_{\theta}H(\theta_t,\gamma_t)\|^2 + \frac{L\alpha_t^2}{2}\|g_{\theta}^t\|^2+\alpha_t \langle \nabla_{\theta}H(\theta_t,\gamma_t), \nabla_{\theta}H(\theta_t,\gamma_t)-g_{\theta}^t \rangle \notag \\
        &\quad -\beta_t\|\nabla_{\gamma}H(\theta_t,\gamma_t)\|^2 + \frac{L\beta_t^2}{2}\|g_{\gamma}^t\|^2+\beta_t \langle \nabla_{\gamma}H(\theta_t,\gamma_t), \nabla_{\gamma}H(\theta_t,\gamma_t)-g_{\gamma}^t \rangle \notag \\ 
        &= H(\theta_t,\gamma_t) -\alpha_t(1-\frac{L\alpha_t}{2})\|\nabla_{\theta}H(\theta_t,\gamma_t)\|^2+\alpha_t(1-L\alpha_t)\langle \nabla_{\theta}H(\theta_t,\gamma_t), \nabla_{\theta}H(\theta_t,\gamma_t)-g_{\theta}^t \rangle
        \notag \\
        &\quad+\frac{L\alpha_t^2}{2}\|g_{\theta}^t-\nabla_{\theta}H(\theta_t,\gamma_t)\|^2 -\beta_t(1-\frac{L\beta_t}{2}) \|\nabla_{\gamma} H(\theta_t,\gamma_t)\|^2 \notag \\
        &\quad+ \beta_t(1-L\beta_t)\langle \nabla_{\gamma}H(\theta_t,\gamma_t), \nabla_{\gamma}H(\theta_t,\gamma_t)-g_{\gamma}^t \rangle +\frac{L\beta_t^2}{2}\|g_{\gamma}^t-\nabla_{\gamma}H(\theta_t,\gamma_t)\|^2. \notag
    \end{align} 
    Let $\delta_{\theta}^t=g_{\theta}^t-\nabla_{\theta}h^{\star}(\theta_t,\gamma_t;\xi_t)$ and $\delta_{\gamma}^t=g_{\gamma}^t-\nabla_{\gamma}h^{\star}(\theta_t,\gamma_t;\xi_t)$ denote the gap caused by $\hat{\xi}^{\star}$ and $\xi^{\star}$, then 
    \begin{align}
         H(\theta_{t+1},\gamma_{t+1}) &\leq H(\theta_t,\gamma_t)-\alpha_t(1-\frac{L\alpha_t}{2})\|\nabla_{\theta}H(\theta_t,\gamma_t)\|^2+\alpha_t(1-L\alpha_t)\langle \nabla_{\theta}H(\theta_t,\gamma_t), \nabla_{\theta}H(\theta_t,\gamma_t)-\nabla_{\theta}h^{\star}(\theta_t,\gamma_t;\xi_t) \rangle \notag \\
         &\quad - \alpha_t(1-L\alpha_t)\langle \nabla_{\theta}H(\theta_t,\gamma_t),\delta_{\theta}^t \rangle+\frac{L\alpha_t^2}{2}\|\nabla_{\theta}h^{\star}(\theta_t,\gamma_t,;\xi_t)+\delta_{\theta}^t-\nabla_{\theta}H(\theta_t,\gamma_t)\|^2 \notag \\
         &\quad -\beta_t(1-\frac{L\beta_t}{2})\|\nabla_{\gamma}H(\theta_t,\gamma_t)\|^2+\beta_t(1-L\beta_t)\langle \nabla_{\gamma}H(\theta_t,\gamma_t), \nabla_{\gamma}H(\theta_t,\gamma_t)-\nabla_{\gamma}h^{\star}(\theta_t,\gamma_t;\xi_t) \rangle
         \notag \\
         &\quad - \beta_t(1-L\beta_t)\langle \nabla_{\gamma}H(\theta_t,\gamma_t),\delta_{\gamma}^t \rangle+\frac{L\beta_t^2}{2}\|\nabla_{\gamma}h^{\star}(\theta_t,\gamma_t,;\xi_t)+\delta_{\gamma}^t-\nabla_{\gamma}H(\theta_t,\gamma_t)\|^2 \notag \\
         &= H(\theta_t,\gamma_t)-\alpha_t(1-\frac{L\alpha_t}{2})\|\nabla_{\theta}H(\theta_t,\gamma_t)\|^2+\alpha_t(1-L\alpha_t)\langle \nabla_{\theta}H(\theta_t,\gamma_t), \nabla_{\theta}H(\theta_t,\gamma_t)-\nabla_{\theta}h^{\star}(\theta_t,\gamma_t;\xi_t) \rangle \notag \\
         &\quad - \alpha_t(1-L\alpha_t)\langle \nabla_{\theta}H(\theta_t,\gamma_t),\delta_{\theta}^t \rangle -\beta_t(1-\frac{L\beta_t}{2})\|\nabla_{\gamma}H(\theta_t,\gamma_t)\|^2  \notag \\
         &\quad +\beta_t(1-L\beta_t)\langle \nabla_{\gamma}H(\theta_t,\gamma_t), \nabla_{\gamma}H(\theta_t,\gamma_t)-\nabla_{\gamma}h^{\star}(\theta_t,\gamma_t;\xi_t) \rangle - \beta_t(1-L\beta_t)\langle \nabla_{\gamma}H(\theta_t,\gamma_t),\delta_{\gamma}^t \rangle\notag \\
         &\quad +\frac{L\alpha_t^2}{2}(\|\delta_{\theta}^t\|^2+\| \nabla_{\theta}h^{\star}(\theta_t,\gamma_t;\xi_t)-\nabla_{\theta}H(\theta_t,\gamma_t)\|^2+2\langle \nabla_{\theta}h^{\star}(\theta_t,\gamma_t;\xi_t)-\nabla_{\theta}H(\theta_t,\gamma_t), \delta_{\theta}^t \rangle) \notag \\
         &\quad +\frac{L\beta_t^2}{2}(\|\delta_{\gamma}^t\|^2+\| \nabla_{\gamma}h^{\star}(\theta_t,\gamma_t;\xi_t)-\nabla_{\gamma}H(\theta_t,\gamma_t)\|^2+2\langle \nabla_{\gamma}h^{\star}(\theta_t,\gamma_t;\xi_t)-\nabla_{\gamma}H(\theta_t,\gamma_t), \delta_{\gamma}^t \rangle). \notag
    \end{align}
    Using $\pm\langle a,b \rangle\leq \frac{1}{2}(\|a\|^2+\|b\|^2)$ in the above inequalities, we can get
    \begin{align}
        H(\theta_{t+1},\gamma_{t+1}) &\leq H(\theta_{t},\gamma_{t}) - \frac{\alpha_t}{2}\|\nabla_{\theta}H(\theta_{t},\gamma_{t})\|^2+\alpha_t(1-L\alpha_t) \langle \nabla_{\theta}H(\theta_t,\gamma_t), \nabla_{\theta}H(\theta_t,\gamma_t) - \nabla_{\theta} h^{\star}(\theta_t,\gamma_t;\xi_t) \rangle \notag \\
        &\quad +\frac{\alpha_t(1+L\alpha_t)}{2}\|\delta_{\theta}^t\|^2 + L\alpha_t^2\|\nabla_{\theta}h^{\star}(\theta_t,\gamma_t;\xi_t)-\nabla_{\theta}H(\theta_t,\gamma_t)\|^2 - \frac{\beta_t}{2}\|\nabla_{\gamma}H(\theta_{t},\gamma_{t})\|^2 \notag \\
        &\quad +\beta_t(1-L\beta_t) \langle \nabla_{\gamma}H(\theta_t,\gamma_t), \nabla_{\gamma}H(\theta_t,\gamma_t) - \nabla_{\gamma} h^{\star}(\theta_t,\gamma_t;\xi_t) \rangle +\frac{\beta_t(1+L\beta_t)}{2}\|\delta_{\gamma}^t\|^2\notag \\
        &\quad + L\beta_t^2\|\nabla_{\gamma}h^{\star}(\theta_t,\gamma_t;\xi_t)-\nabla_{\gamma}H(\theta_t,\gamma_t)\|^2.\notag
    \end{align}
    In each iteration of proposed algorithm, $\hat{\xi^{\star}_t}$ is generated to calculate the gradient, then error $\delta_{\theta}^t$ and $\delta_{\gamma}^t$ satisfy
    \begin{align}
        \|\delta_{\theta}^t\|^2 \leq \|\nabla_{\theta}h(\theta_t,\gamma_t;\xi_t,\xi^{\star}_t)-\nabla_{\theta}h(\theta_t,\gamma_t;\xi_t,\hat{\xi_t^{\star}})\|^2 \leq L_{\theta\xi}^2\|\xi^{\star}_t-\hat{\xi^{\star}_t}\|^2 \leq \frac{2L_{\theta\xi}^2}{\mu_h}\varepsilon, \notag
    \end{align}
    \begin{align}
        \|\delta_{\gamma}^t\|^2 \leq \|d(\xi_t,\xi^{\star}_t)-d(\xi_t,\xi^{\star}_t)\|^2 \leq \|\xi^{\star}_t-\hat{\xi^{\star}_t}\|^2 \leq \frac{2}{\mu_h}\varepsilon. \notag
    \end{align}
    Because of exchangeable principle between gradient and expectation, we have $\mathbb{E}_Q[\nabla H-\nabla h^{\star}]=0$. Then the inequality can be simplified as
    \begin{align}
        \mathbb{E}[H(\theta_{t+1},\gamma_{t+1})-H(\theta_t,\gamma_t)] &\leq - \frac{\alpha_t}{2}\|\nabla_{\theta}H(\theta_t,\gamma_t)\|^2+\frac{\alpha_t(1+L\alpha_t)L_{\theta\xi}^2}{\mu_h}\varepsilon+L\alpha_t^2\sigma_{\theta}^2 \notag \\
        &\quad -\frac{\beta_t}{2}\|\nabla_{\gamma}H(\theta_t,\gamma_t)\|^2+\frac{\beta_t(1+L\beta_t)}{\mu_h}\varepsilon+L\beta_t^2\sigma_{\gamma}^2. \notag
    \end{align}
    Considering a fixed stepsize $\alpha_t=\beta_t=\alpha$, and calculating the total gap in $T$ iterations training, we have
    \begin{align}
        \frac{1}{T}\sum^{T}_{t=0}\mathbb{E}_{Q}[ \|\nabla_{\theta}H(\theta_t,\gamma_t)\|^2 + \|\nabla_{\gamma}H(\theta_t,\gamma_t)\|^2 &\leq 
        \frac{2}{\alpha T}[H(\theta_0,\gamma_0)-H(\theta_T,\gamma_T)] + 2L\alpha(\sigma_{\theta}^2 + \sigma_{\gamma}^2) + \frac{2(1+L\alpha)}{\mu_h}(L_{\theta\xi}^2+1)\varepsilon \notag \\
        &\quad\leq \frac{2}{\alpha T}\Lambda + 2L\alpha(\sigma_{\theta}^2 + \sigma_{\gamma}^2) + \frac{2(1+L\alpha)}{\mu_h}(L_{\theta\xi}^2+1)\varepsilon. \notag
  \end{align}
  Let $\alpha=\sqrt{\frac{\Lambda}{L(\sigma_{\theta}^2+\sigma_{\gamma}^2)T}}$ and $\varepsilon=\frac{\mu_h}{(L_{\theta\xi}^2+1)\sqrt{T}}$, then the bound should be
  \begin{align}
       \frac{1}{T}\sum^{T}_{t=0}\mathbb{E}_{Q}[ \|\nabla_{\theta}H(\theta_t,\gamma_t)\|^2 + \|\nabla_{\gamma}H(\theta_t,\gamma_t)\|^2 ] \leq 
        (4\sqrt{L(\sigma_{\theta}^2+\sigma_{\gamma}^2)\Lambda} +2)\sqrt{\frac{1}{T}} +
        2\sqrt{\frac{L\Lambda}{\sigma_{\theta}^2+\sigma_{\gamma}^2}}\frac{1}{T}. \notag
  \end{align}
  And this gives the result.
    \subsection{Proof of Corollary \ref{col2}} \label{secA.4}
    If $c(f(\theta;X);Y)$ is $\lambda-$strongly convex in $\theta$, then
    \begin{align}
        h(\theta,\gamma;\xi,\xi')-h(\theta',\gamma;\xi,\xi') \geq \langle \nabla_{\theta}h(\theta'
        ,\gamma;\xi,\xi'), \theta-\theta' \rangle+ \frac{\lambda}{2}\|\theta-\theta'\|^2.  \notag
    \end{align}
    And $h(\theta,\gamma;\xi,\xi')$ is linear in $\gamma$, then we have
    \begin{align}
        h(\theta_2,\gamma_2;\xi,\xi')-h(\theta_1,\gamma_1;\xi,\xi') &= h(\theta_2,\gamma_2;\xi,\xi')-h(\theta_2,\gamma_1;\xi,\xi') + h(\theta_2,\gamma_1;\xi,\xi')-h(\theta_1,\gamma_1;\xi,\xi') \notag \\
        &\geq \langle \nabla_{\gamma}h(\theta_2,\gamma_1;\xi,\xi'), \gamma_2-\gamma_1 \rangle + \langle \nabla_{\theta}h(\theta_1,\gamma_1;\xi,\xi'), \theta_2-\theta_1 \rangle + \frac{\lambda}{2}\|\theta_2-\theta_1\|^2 \notag \\
        &=\langle \nabla_{\gamma}h(\theta_1,\gamma_1;\xi,\xi'), \gamma_2-\gamma_1 \rangle + \langle \nabla_{\theta}h(\theta_1,\gamma_1;\xi,\xi'), \theta_2-\theta_1 \rangle + \frac{\lambda}{2}\|\theta_2-\theta_1\|^2. \notag
    \end{align}
    Then based on the definition of $\xi^{\star}$, the equality becomes
    \begin{align}
        h(\theta_2,\gamma_2;\xi,\xi^{\star}_2) &\geq 
        h(\theta_2,\gamma_2;\xi,\xi^{\star}_1)  \notag \\
        &\geq
        h(\theta_1,\gamma_1;\xi,\xi^{\star}_1) + \langle \nabla_{\gamma}h(\theta_1,\gamma_1;\xi,\xi'), \gamma_2-\gamma_1 \rangle + \langle \nabla_{\theta}h(\theta_1,\gamma_1;\xi,\xi'), \theta_2-\theta_1 \rangle + \frac{\lambda}{2}\|\theta_2-\theta_1\|^2. \notag
    \end{align}
    Thus, the following two inequalities can be generated by this
    \begin{align}
        h(\theta^{\star},\gamma^{\star};\xi,\xi^{\star}_t) \geq
        h(\theta_t,\gamma_t;\xi,\xi^{\star}_t) + \langle \nabla_{\gamma}h(\theta_t,\gamma_t;\xi,\xi^{\star}_t), \gamma^{\star}-\gamma_t \rangle + \langle \nabla_{\theta}h(\theta_t,\gamma_t;\xi,\xi^{\star}_t), \theta^{\star}-\theta_t \rangle + \frac{\lambda}{2}\|\theta^{\star}-\theta_t\|^2, \notag
    \end{align}
    \begin{align}
        h(\theta_t,\gamma_t;\xi,\xi^{\star}_t) \geq
        h(\theta^{\star},\gamma^{\star};\xi,\xi^{\star}_t) - \langle \nabla_{\gamma}h(\theta^{\star},\gamma^{\star};\xi,\xi^{\star}_t), \gamma^{\star}-\gamma_t \rangle+ \frac{\lambda}{2}\|\theta^{\star}-\theta_t\|^2. \notag
    \end{align}
    After getting these inequalities, the distance between solutions obtained from algorithm and optimal solutions should be
    \begin{align}
        \mathbb{E}_{Q}[\|\theta_{t+1}-\theta^{\star}\|^2 + \|\gamma_{t+1}-\gamma^{\star}\|^2] 
        &= \mathbb{E}_{Q}[\|\theta_{t}-\alpha g_{\theta}^t-\theta^{\star}\|^2 + \|\gamma_{t}-\alpha g_{\gamma}^t-\gamma^{\star}\|^2] \notag \\
        &= \mathbb{E}_{Q}[\|\theta_{t}-\theta^{\star}\|^2 + \|\gamma_{t}-\gamma^{\star}\|^2] -2\alpha \mathbb{E}_Q[ \langle g_{\theta}^t, \theta_t-\theta^{\star} \rangle + \langle g_{\gamma}^t, \gamma_t-\gamma^{\star} \rangle] \notag \\
        &\quad + \alpha^2 \mathbb{E}_Q[\|g_{\theta}^t\|^2+\|g_{\gamma}^t\|^2] \notag \\
        &\leq \mathbb{E}_{Q}[\|\theta_{t}-\theta^{\star}\|^2 + \|\gamma_{t}-\gamma^{\star}\|^2]+ \alpha^2 \mathbb{E}_Q[\|g_{\theta}^t\|^2+\|g_{\gamma}^t\|^2] \notag \\ &\quad -2\alpha \mathbb{E}_Q[ h(\theta_t,\gamma_t;\xi,\xi^{\star}_t) - h(\theta^{\star},\gamma^{\star};\xi,\xi^{\star}_t) + \frac{\lambda}{2}\|\theta_t-\theta^{\star}\|^2] \notag \\
        &\leq \mathbb{E}_{Q}[\|\theta_{t}-\theta^{\star}\|^2 + \|\gamma_{t}-\gamma^{\star}\|^2]+ \alpha^2 \mathbb{E}_Q[\|g_{\theta}^t\|^2+\|g_{\gamma}^t\|^2]-\alpha \lambda \mathbb{E}_Q[\|\theta_t-\theta^{\star}\|^2]. \notag
    \end{align}
    The final inequality holds because of the definition of $\theta^{\star}$ and $\gamma^{\star}$. $\theta$ is more important than $\gamma$, so we only consider the gap between $\theta_t$ and $\theta^{\star}$.
    \begin{align}
        \mathbb{E}_Q[\|\theta_{t+1}-\theta^{\star}\|^2] 
        &\leq \mathbb{E}_Q[\|\gamma_t-\gamma^{\star}\|^2-\|\gamma_{t+1}-\gamma^{\star}\|^2] + (1-\alpha\lambda) \mathbb{E}_Q [\|\theta_t-\theta^{\star}\|^2] + \alpha^2 \mathbb{E}_Q[\|g_{\theta}^t\|^2+\|g_{\gamma}^t\|^2] \notag \\
        &\leq \mathbb{E}_Q[\|\gamma_t-\gamma_{t+1}\|^2]+ (1-\alpha\lambda) \mathbb{E}_Q [\|\theta_t-\theta^{\star}\|^2] + \alpha^2 \mathbb{E}_Q[\|g_{\theta}^t\|^2+\|g_{\gamma}^t\|^2] \notag \\
        &\leq (1-\alpha\lambda) \mathbb{E}_Q [\|\theta_t-\theta^{\star}\|^2] + 2\alpha^2 \mathbb{E}_Q[\|g_{\theta}^t\|^2+\|g_{\gamma}^t\|^2]. \notag
    \end{align}
    To simplify the inequality, let $F_{t}=\mathbb{E}_Q[\|\theta_{t}-\theta^{\star}\|^2]$, then 
    \begin{equation}
        F_{t+1} \leq (1-\alpha\lambda)F_t + 2\alpha^2\mathbb{E}_Q[\|g_{\theta}^t\|^2+\|g_{\gamma}^t\|^2]. \notag
    \end{equation}
    From Theorem \ref{thm1}, we have known that $\alpha=\frac{A}{\sqrt{T}}$, where $A$ is a constant whose value described in Theorem \ref{thm1}, and from Assumption \ref{asp3}, we can get
    \begin{align}
        \mathbb{E}_Q[\|g_{\theta}^t\|^2+\|g_{\gamma}^t\|^2] \leq \mathbb{E}_Q [{\nabla_{\theta}H(\theta_t,\gamma_t)+ \nabla_{\gamma}H(\theta_t,\gamma_t)}] +[ \sigma_{\theta}^2 + \sigma_{\gamma}^2] \leq M, \notag
    \end{align}
    where $M$ is a large value, when $\xi'_t=\xi^{\star}_t$ in $h(\theta_t,\gamma_t;\xi_t,\xi'_t)$. Then the recursive formula should be
    \begin{align}
        F_{t+1} \leq (1-\frac{A}{\sqrt{T}})F_t + \frac{M'}{T}, \notag
    \end{align}
    where $M'=MA^2$ is a constant. When $A\geq1$, $A$ can be replaced by 1 because of inequality, and when $A<1$, if $T$ is large enough, $1-\frac{A}{\sqrt{T}}\approx 1-\frac{1}{\sqrt{T}}$. Thus,  
    \begin{align}
        F_{t+1} 
        &\leq (1-\frac{1}{\sqrt{T}})F_t+\frac{M'}{T} \notag \\
        &\leq (1-\frac{1}{\sqrt{T}})^2F_{t-1}+ \frac{M'}{T}(1+(1-\frac{1}{\sqrt{T}})) \notag \\
        &\leq (1-\frac{1}{\sqrt{T}})^{T+1}F_0+ \frac{M'}{T} \sum_{i=0}^{T-1}(1-\frac{1}{\sqrt{T}})^i, \notag \\
        &\leq \frac{1}{\sqrt{T}}F_0 + \frac{M'}{\sqrt{T}}. \notag
    \end{align}
    To prove the final inequality, we need to show $(1-\frac{1}{\sqrt{T}})^{T}\leq \frac{1}{\sqrt{T}}$. 
    \begin{align}
        (1-\frac{1}{\sqrt{T}})^{T}&\leq \frac{1}{\sqrt{T}} \notag \\
        T\ln{(1-\frac{1}{\sqrt{T}})} &\leq \ln{\frac{1}{\sqrt{T}}}, \notag
    \end{align}
    Let $u=1/\sqrt{T}$, then the inequality becomes
    \begin{align}
        u^2\ln u - \ln(1-u) \geq 0. \notag
    \end{align}
    Let $g(u)=u^2\ln u - \ln(1-u)$, and $\lim_{u\rightarrow 0}g(u)=0$. So $g(u)\geq 0 \Leftrightarrow g'(u)\geq 0$.
    \begin{align}
        2u\ln u + u &+ \frac{1}{1-u} \geq 0 \notag \\
        &\Updownarrow \notag \\
        2\ln u + 1 +& \frac{1}{u} + \frac{1}{1-u} \geq 0, \notag
    \end{align}
    which is easy to prove.
    Therefore, $\mathbb{E}_Q[\|\theta_{T}-\theta^{\star}\|^2]$ is bounded by $O(1/\sqrt{T})$.
    
    \section{Proof of Section \ref{sec4.2}} \label{secB}
    \subsection{Proof of Theorem \ref{probguarantee}} \label{secB.1}
    For any distribution $P_1$ and $P_2$, Wasserstein distance with norm-form mass transportation cost can be controlled by weighted total variation \citep[Theorem 6.15]{villani2009optimal}, which is
    \begin{equation}
        W_d(P_1,P_2) \leq \int_{\Xi^2} d(\xi_1,\xi_2) d|P_1-P_2|(\xi) = \|d(\xi_1,\xi_2)(P_1-P_2)\|_{TV}. \notag
    \end{equation}
    If Assumption \ref{asp1} holds, this inequality can be relaxed by diameter $D$ of $\Xi$,
    \begin{equation}
        W_d(P_1,P_2) \leq \int_{\Xi^2} D d|P_1-P_2|(\xi) = D\|P_1-P_2\|_{TV}. \notag
    \end{equation}
    Then the weighted total variation (WTV) is relaxed as a total variation (TV). Based on weighted CKP inequalities \citep[Theorem 2.1 and Remark 2.2]{bolley2005weighted}, we get three inequalities about Wasserstein distance
    \begin{align} \label{CKP1}
        \|d(\xi_1,\xi_2)(P_1-P_2)\|_{TV} &\leq \Big(\frac{3}{2} + \log \int e^{2d(\xi_1,\xi_2)} dP_2(\xi_2) \Big) \Big( \sqrt{H(P_1|P_2)}+\frac{1}{2}H(P_1|P_2) \Big) \notag \\
        &\leq \Big(\frac{3}{2}+2D\Big) \Big(\sqrt{H(P_1|P_2)} +\frac{1}{2}H(P_1|P_2) \Big),
    \end{align}
    \begin{align} \label{CKP2}
        \|d(\xi_1,\xi_2)(P_1-P_2)\|_{TV} &\leq \sqrt{2} \Big( 1+\log \int e^{d(\xi_1,\xi_2)^2} dP_2(\xi_2) \Big)^{1/2} \sqrt{H(P_1|P_2)} \notag \\
        &\leq \sqrt{2} \Big(1+D^2 \Big)^{1/2}\sqrt{H(P_1|P_2)},
    \end{align}
    \begin{align} \label{CKP3}
        \|P_1-P_2\|_{Tv} \leq 2\sqrt{H(P_1|P_2)}.
    \end{align}
    And we have 
    \begin{equation}
        \|d(\xi_1,\xi_2)(P_1-P_2)\|_{TV} \leq D\|P_1-P_2\|_{TV}. \notag
    \end{equation}
    If we consider to construct the complement of the ambiguity set, that is selecting a radius $\rho$ such that $\rho<W_d(P_1,P_2)$, then
    \begin{equation} 
        \rho<\|P_1-P_2\|_{TV}\leq \left\{
        \begin{aligned}
        &\Big(\frac{3}{2}+2D\Big) \Big (\sqrt{H(P_1|P_2)} +\frac{1}{2}H(P_1|P_2)\Big),\\
        &\sqrt{2}(1+D^2)^{1/2}\sqrt{H(P_1|P_2)},\\
        &2D\sqrt{H(P_1|P_2)}.\\
        \end{aligned} \notag
        \right.
    \end{equation}    
    By solving these inequalities, we can get the lower bound of $H(P_1|P_2)$
    \begin{equation}
    H(P_1|P_2)\geq
    \left\{
    \begin{aligned}
    &\Big( \sqrt{1+\frac{4\rho}{4D+3}}-1 \Big)^2 \quad &&if\; &&D\geq 1 \ and \  \rho\geq 8D-\frac{12D}{4D+3} \ or\\ 
    & && &&D\leq1\ and\ \rho\geq \frac{8(1+D^2)}{4D+3}+2\sqrt{2(1+D^2)},\\
    &\frac{\rho^2}{2(1+D^2)}\quad &&if\; &&D\geq 1\ and \ \rho<8D-\frac{12D}{4D+3},\\
    &\frac{\rho^2}{4D^2}&&if &&D\leq 1\ and \ \rho<\frac{8(1+D^2)}{4D+3}+2\sqrt{2(1+D^2)}.\\
    \end{aligned} \notag
    \right.
    \end{equation}
    When $\rho<D$, the condition of the first inequality cannot be satisfied. Therefore, the lower bound becomes
    \begin{equation}
    H(P_1|P_2)\geq
    \left\{
    \begin{aligned}
    &\frac{\rho^2}{2(1+D^2)}\quad &&if\; &&D\geq 1 ,\\
    &\frac{\rho^2}{4D^2}&&if &&D < 1 .\\
    \end{aligned} \notag
    \right.
    \end{equation}
    \cite{ji2021data} combine $H(P_1|P_2)$ with Wasserstein distance to derive the conclusion
    \begin{equation}
        P(W_d(P_{true},P_N)\leq \rho_N) \geq 1- e^{-N \inf H(P_N|P)}. \notag
    \end{equation}
    Then we can get the results
    \begin{equation} \notag
		P(W_d(P_{true},P_N)\leq\rho_N)\geq
		\left\{
		\begin{aligned}
			&1-e^{-N\frac{\rho_N^2}{2(1+D^2)}}\quad &&if\; &&D\geq 1,\\
			&1-e^{-N\frac{\rho_N^2}{4D^2}}&&if &&D < 1.\\
		\end{aligned}
		\right.
	\end{equation}
    \subsection{Proof of Theorem \ref{finiesampleguarantee}} \label{secB.2}
    $r_{DRO}$ can be calculated by 
    \begin{align}
    r_{DRO} &= [C(\theta_R)-C_{DRO}(\theta_R;\rho_N)]+[C_{DRO}(\theta_R;\rho_N)-C_{DRO}(\theta^{\star};\rho_N)] \notag\\
        &\ \ +[C_{DRO}(\theta^{\star};\rho_N)-C_N(\theta^{\star})]+[C_N(\theta^{\star})-C(\theta^{\star})], \notag
    \end{align}
    and it is divided into four parts. In \cite{gao2023finite}, parts of the results have been proved. In Section \ref{secB.1}, Wasserstein distance have been bounded by $H(P_{true}|P_N)$, with $\tau=4D^2$ when $D<1$ and $\tau=2(1+D^2)$ when $D\geq 1$. Then we have the following result \citep[Corollary 1]{gao2023finite} with probability at least $1-e^{-t}$
    \begin{equation}
        C(\theta) \leq  C_N(\theta) + \sqrt{\frac{\tau t}{N}}L_c \leq C_{DRO}(\theta) + \sqrt{\frac{\tau t}{N}}L_c, \notag
    \end{equation}
    where the value of $\tau$ is defined before, $L_c$ is the Lipschitz constant of $c(f(\theta;X);Y)$ in $\xi$. And the first part is bounded by $O(1/\sqrt{N})$.

    Based on the definition of $\theta_R$, it is obvious that 
    \begin{equation}
        C_{DRO}(\theta_R;\rho_N) \leq C_{DRO} (\theta^{\star;\rho_N}), \notag
    \end{equation}
    which means the second part is bounded by 0. The third one can be bounded by Lipschitz regularization method \citep[Lemma 1]{gao2023finite},
    \begin{equation}
        C_{DRO}(\theta;\rho) - C_N{DRO}(\theta) \leq \rho L_c.
    \end{equation}
    Then after applied the radius $\rho_N$ chosen in Corollary \ref{radiuschosen} with probability $q$, its bound is also $O(1/\sqrt{N})$. 

    The final part is the gap between empirical distribution and ground-truth distribution, which is also known as the performance of sample average approximation (SAA). There have been many researches about this gap, which satisfies $C_N(\theta^{\star})-C(\theta^{\star})=O_P(1/\sqrt{N})$ under mild conditions \citep[Theorem 8.3]{fu2015handbook}. Thus the gap $r_{DRO}$ is bounded by $O(1/\sqrt{N})$ with probability $q(1-e^{-t})$.
    
    \section{Proof of Section \ref{sec4.3}}
    In \cite{shalev2012online} Lemma 2.1, for arbitrary function $f(x)$, we have known 
    \begin{align}
        \sum_{t=1}^T(f(x_t)-f(u)) \leq \sum_{t=1}^T(f(x_t)-f(x_{t+1})). \notag
    \end{align}
    \subsection{Proof of Corollary \ref{regret}}
    Based on the definition of regret, we have
    \begin{align}
        \mathbb{E}_Q [Regret_T(\theta^{\star},\gamma^{\star})] 
        &= \mathbb{E}_Q \Big[ \sum_{t=1}^T \big(h(\theta_t,\gamma_t;\xi_t,\xi_t^{\star})-h(\theta^{\star},\gamma^{\star};\xi_t,\xi^{\star}) \big) \Big] \notag \\
        &\leq \mathbb{E}_Q \Big[ \sum_{t=1}^T \big(h(\theta_t,\gamma_t;\xi_t,\xi_t^{\star})-h(\theta_{t+1},\gamma_{t+1};\xi_t,\xi_{t+1}^{\star}) \big) \Big] \notag \\
        &\leq \mathbb{E}_Q \Big[ \sum_{t=1}^T \big(h(\theta_t,\gamma_t;\xi_t,\xi_t^{\star})-h(\theta_{t+1},\gamma_{t+1};\xi_t,\xi_t^{\star}) \big) \Big] \notag \\
        &\leq \mathbb{E}_Q \Big[ \sum_{t=1}^T \big[ \big(h(\theta_t,\gamma_t;\xi_t,\xi_t^{\star})-h(\theta_{t+1},\gamma_t;\xi_t,\xi_t^{\star}) \big) + \big(h(\theta_{t+1},\gamma_t;\xi_t,\xi_t^{\star})-h(\theta_{t+1},\gamma_{t+1};\xi_t,\xi_t^{\star}) \big) \big]  \Big] \notag \\
        &\leq \mathbb{E}_Q \Big[ \sum_{t=1}^T \alpha_t\big( \|g_{\theta}^t\|^2 + \|g_{\gamma}^t\|^2 \big)  \Big] \notag \\
        &\leq \mathbb{E}_Q \Big[ \sum_{t=1}^T \alpha_t\big( (\nabla_{\theta}H(\theta_t,\gamma_t)+ \nabla_{\gamma}H(\theta_t,\gamma_t))+(\sigma_{\theta}^2+\sigma_{\gamma}^2)\big)  \Big], \notag
    \end{align}
    where the second inequality holds because of the definition of $\xi^{\star}$, the forth inequality holds because of the convexity and the final one holds because of Assumption \ref{asp4}.In Theorem \ref{thm1}, $\alpha_T$ is set as a constant step size and satisfies $\alpha \propto \frac{1}{\sqrt{T}}$. Thus, we have 
    \begin{equation}
         \mathbb{E}_Q [Regret_T(\theta^{\star},\gamma^{\star})] \leq O(\sqrt{T}) \notag.
    \end{equation}
    And the convergence speed is $O(\frac{1}{\sqrt{T}})$.
\end{APPENDICES}

\end{document}